\documentclass[onefignum,onetabnum]{siamart220329}

\usepackage{amssymb}
\usepackage{graphicx}
\usepackage{booktabs}
\usepackage{needspace}
\usepackage{enumitem}

\newsiamthm{lem}{Lemma}
\newsiamthm{prop}{Proposition}
\newsiamremark{assum}{Assumption}
\newsiamremark{algo}{Algorithm}
\newsiamremark{rem}{Remark}

\newcommand{\bsigma}{\mbox{\boldmath$\sigma$}}
\newcommand{\btau}{\mbox{\boldmath$\tau$}}
\newcommand{\bchi}{\mbox{\boldmath$\chi$}}
\newcommand{\beps}{\mbox{\boldmath$\varepsilon$}}
\newcommand{\bb}{\mbox{\boldmath$b$}}
\newcommand{\bSigma}{{\mathbf{\Sigma}}}
\newcommand{\bff}{{\mathbf{f}}}
\newcommand{\bn}{{\mathbf{n}}}
\newcommand{\bzero}{{\mathbf{0}}}
\newcommand{\cA}{{\mathcal{A}}}
\newcommand{\cD}{{\mathcal{D}}}
\newcommand{\cJ}{{\mathcal{J}}}
\newcommand{\cR}{{\mathcal{R}}}
\newcommand{\cT}{{\mathcal{T}}}
\newcommand{\grad}{\nabla}
\newcommand{\divv}{\nabla\cdot}
\DeclareMathOperator{\Span}{Span}
\newcommand{\divvr}{\operatorname{div}}
\newcommand{\tri}{|\!|\!|}
\newcommand{\Hdiv}{H(\divvr;\,\Omega)}
\newcommand{\R}{\mathbb{R}}
\newcommand{\cX}{{\mathcal{X}}}
\newcommand{\cY}{{\mathcal{Y}}}
\newcommand{\norm}[1]{\left\lVert #1\right\rVert}
\newcommand{\tnorm}[1]{\tri #1\tri}
\DeclareMathOperator{\tr}{tr}
\DeclareMathOperator{\dev}{dev}
\DeclareMathOperator{\skw}{skw}
\newcommand{\bI}{{\mathbf{I}}}
\newcommand{\bM}{{\mathbf{M}}}
\newcommand{\bR}{{\mathbf{R}}}
\newcommand{\bu}{{\mathbf{u}}}
\newcommand{\bv}{{\mathbf{v}}}
\newcommand{\aLS}{\alpha_{\rm LS}}
\newcommand{\bLS}{\beta_{\rm LS}}
\newcommand{\rb}{{\rm rb}}
\newcommand{\fe}{{\rm fe}}
\newcommand{\snap}{{\rm snap}}
\newcommand{\train}{{\rm train}}
\newcommand{\LS}{{\rm LS}}
\newcommand{\Hmoo}{H^{-1/2}_{00}}
\newcommand{\rrf}{\mathtt{r}_{\rb,\fe}}
\newcommand{\bbeta}{\mbox{\boldmath$\beta$}}

\makeatletter
\AtBeginDocument{%
  \def\refstepcounter@optarg[#1]#2{%
    \cref@old@refstepcounter{#2}%
    \cref@constructprefix{#2}{\cref@result}%
    \@ifundefined{cref@#1@alias}{\def\@tempa{#1}}%
      {\def\@tempa{\csname cref@#1@alias\endcsname}}%
    \protected@edef\cref@currentlabel{[\@tempa][\arabic{#2}][\cref@result]%
      \csname p@#2\endcsname\csname the#2\endcsname}}}
\makeatother

\title{Non-Intrusive Least-Squares Reduced Basis Methods: Black-Box
Snapshots and Exact Residual Certificates\thanks{\funding{This work was supported in part by the Research Grants Council of
the Hong Kong SAR, China, under the GRF Grant Projects No.\ CityU 11316222
and CityU 11305025.}}}
\author{Shun Zhang\thanks{Department of Mathematics, City University of Hong
Kong, Kowloon Tong, Hong Kong, China (\email{shun.zhang@cityu.edu.hk}).}}
\headers{Non-intrusive least-squares reduced basis methods}{Shun Zhang}

\begin{document}
\maketitle

\begin{abstract}
For a parametrized linear first-order system, the least-squares functional
at a conforming pair is the squared norm of the residual of the continuous
problem, and it can be computed exactly.  In this paper, we use this
least-squares residual as the error certificate of reduced basis methods whose snapshots
are computed by black-box solvers.  The output of the solver is first made
conforming by a reconstruction, and a flux is recovered by least-squares if
the solver does not provide one.  The online reduced problem minimizes the least-squares functional over the product of the
space spanned by the snapshot fluxes and the space spanned by the snapshot
potentials.  It needs no inf--sup condition on the
reduced spaces, and its certificate needs no stability constant.  We discuss second-order elliptic equations, the Stokes equations,
where only the velocity of a black-box solver is used, and linear
transport equations with a moving inflow boundary.  In the numerical
experiments, the snapshots are computed by least-squares, conforming,
nonconforming, stabilized, finite volume, and finite difference solvers.  With adaptive meshes
and the same tolerances, the greedy algorithm stops at the same or nearly
the same reduced dimension for all solvers that accept a mesh.
\end{abstract}

\begin{keywords}
reduced basis method, least-squares, non-intrusive, a posteriori error
estimate, adaptive greedy algorithm
\end{keywords}

\begin{MSCcodes}
65N30, 65N15, 65N50, 76D07
\end{MSCcodes}

\section{Introduction}

Reduced basis (RB) methods approximate the solutions of a parametrized
partial differential equation for many parameter values in the span of a
few precomputed solutions, called snapshots; see the review articles and
books~\cite{BennerEtAl2017,Haasdonk2017,HesthavenRozzaStamm2016,Maday2006,OhlbergerRave2016,QuarteroniManzoniNegri2016,QuarteroniRozzaManzoni2011,RozzaHuynhPatera2008}
and the recent
surveys~\cite{HesthavenPagliantiniRozza2022,HesthavenPeherstorferUnger2026}
on time-dependent and transport-dominated problems.  An RB method has an
offline stage, where the snapshots are computed, and an online stage, where
a small problem on the span of the snapshots is solved for each new
parameter.  We call this small problem the online reduced problem.  Its cost
does not depend on the dimension of the high-fidelity discretization.  An a~posteriori error estimate certifies the RB
approximation and drives the greedy selection of the snapshots.  In its
classical form, the estimate bounds the error with respect to a discrete
truth solution by a dual norm of the residual divided by a lower bound of a
coercivity or inf--sup
constant~\cite{HesthavenRozzaStamm2016,RozzaHuynhPatera2008}.  In this
paper, the error estimate is the norm of the residual of the continuous
first-order system, that is, the square root of a least-squares functional.
It is computed exactly, without a truth discretization, and no stability
constant is needed.

\subsection{Two properties of the least-squares functional}

Let $\cD\subset\R^p$ be a compact parameter domain.  For $\mu\in\cD$, we
write the parametric problem as a first-order system $\cA_\mu U(\mu)=F_\mu$.
Here $\cX$ is a real Hilbert space of pairs $W=(\btau,v)$ of a flux $\btau$
and a potential $v$ on a domain $\Omega\subset\R^d$, for example,
$\Hdiv\times H^1(\Omega)$ with homogeneous
essential boundary conditions; $\cY$ is a product of $L^2$ spaces;
$\cA_\mu:\cX\to\cY$ is linear, bounded, and injective; $F_\mu\in\cY$; and
$U(\mu)=(\bsigma(\mu),u(\mu))\in\cX$ is the exact solution.  Nonhomogeneous
essential boundary data are lifted (section~\ref{sec:abstract}), and the
pairs in $\cX$, shifted by the lifting, are called conforming pairs.  For
$W\in\cX$, define
\begin{equation}
  \tnorm{W}_\mu:=\norm{\cA_\mu W}_{\cY}
  \quad\text{and}\quad
  \eta(W;\mu):=\norm{\cA_\mu W-F_\mu}_{\cY}.
\end{equation}
Since $F_\mu=\cA_\mu U(\mu)$, we have
\begin{equation}\label{eq:exactcert}
  \eta(W;\mu)=\tnorm{U(\mu)-W}_\mu
  \qquad\text{for all }W\in\cX .
\end{equation}
This identity gives the following two properties.
\begin{enumerate}[label=(P\arabic*)]
\item\label{it:P1} \emph{As a method.}  On any finite-dimensional subspace
$X\subset\cX$, the minimizer of $\eta(\cdot;\mu)^2$ exists, is unique, and
is the orthogonal projection of $U(\mu)$ onto $X$ with respect to the inner
product $(\cA_\mu\cdot,\cA_\mu\cdot)_{\cY}$.  No inf--sup condition on $X$
is needed.
\item\label{it:P2} \emph{As a certificate.}  For every $W\in\cX$,
$\eta(W;\mu)$ is the norm of the residual of the continuous problem at $W$,
and by~\eqref{eq:exactcert}, it is the error of $W$ with respect to the
exact solution in $\tnorm{\cdot}_\mu$.  The function $W$ does not need
to be the least-squares solution or the solution of any discrete problem.
\end{enumerate}
Both properties are classical for least-squares finite element
methods~\cite{BochevGunzburger2009,CaiLazarovManteuffelMcCormick1994,Zhang2023}.
Property~\ref{it:P2} is the basis of the non-intrusive least-squares
functional error estimator.  This estimator was proposed
in~\cite{CaiZhang2010} for conforming Galerkin solutions, where the flux is
recovered by minimizing the functional with the Galerkin solution fixed, and
was extended in~\cite{LiZhang2025} to general second-order elliptic and
nonlinear problems; see~\cite{ZhangWeak2026} for weakly imposed boundary
conditions.

\subsection{The online and the offline stages}

In our method, least-squares is used in the online stage, where both
properties are needed.  The online reduced problem minimizes the functional
over the reduced space; by~\ref{it:P1}, no inf--sup condition is needed, and
by~\ref{it:P2}, the square root of the minimum certifies the reduced
solution.  The snapshots can be computed
by any solver whose output can be made conforming, since~\ref{it:P2} holds
for every conforming pair.  For a selected parameter $\mu$, the offline
computation of one snapshot is
\begin{equation}
  u_h(\mu)
  \ \overset{\cR_h}{\longrightarrow}\
  \widetilde u_h(\mu)
  \ \overset{\text{flux recovery}}{\longrightarrow}\
  \bsigma_h(\mu)
  \ \longrightarrow\
  \bigl(\bsigma_h(\mu),\widetilde u_h(\mu)\bigr) .
\end{equation}
Here $u_h(\mu)$ is computed by the user's solver, $\cR_h$ is a conforming
reconstruction (the identity for a conforming solver), and $\bsigma_h(\mu)$
minimizes the functional with $\widetilde u_h(\mu)$ fixed, unless the solver
provides a conforming flux.  We introduce four levels of
non-intrusiveness; see Table~\ref{tab:levels}.  The user
implements the least-squares forms and the flux recovery, and the black-box
solver is called on a given mesh or, at non-intrusive level~3, on its own
grid.  The
truth solver only determines the snapshots; the reduction and the
certification are done by least-squares.  In model order reduction, the
term ``non-intrusive'' usually refers to a regression of the reduced
coefficients on the
parameter~\cite{HesthavenUbbiali2018,KastGuoHesthaven2020,PadulaGirfoglioRozza2024}.
Here, only the snapshots are computed non-intrusively, and the
online stage is a projection with an exact certificate.

\begin{table}[htbp]
\centering
\small
\begin{tabular}{llll}
\toprule
non-intrusive level & truth solver & mesh & $\cR_h$ \\
\midrule
0 & least-squares                    & given        & identity   \\
1 & conforming Galerkin              & given        & identity   \\
2 & DG, nonconforming, finite volume & given        & needed     \\
3 & finite difference, legacy code   & solver's own & needed     \\
\bottomrule
\end{tabular}
\caption{Non-intrusive levels for a potential in $H^1$ (DG:
discontinuous Galerkin).  The online reduced problem and the certificate
are the same for all four levels; only $\cR_h$ changes.}
\label{tab:levels}
\end{table}

\subsection{Related work}

In~\cite{Yano2015,Yano2016,Yano2018}, Yano
developed RB methods that certify the error with respect to the exact
solution and adapt the truth discretization.  For symmetric coercive
equations, upper and lower bounds of the energy of the exact solution
follow from the primal and the complementary variational
principles~\cite{Yano2015}.  The minimum-residual mixed RB method
of~\cite{Yano2016} is a least-squares RB method: the snapshots and the
reduced solution minimize the least-squares functional of the first-order
system, the primal and flux reduced spaces are decoupled, and the snapshots
are kept on the refined meshes.  The square root of the functional bounds the
dual norm of the residual of the continuous problem, and a lower bound of the
inf--sup constant turns it into a bound of the error.  Energy-norm
and output error bounds with spatio-parameter adaptivity are derived
in~\cite{Yano2018}.  The primal--dual RB
method~\cite{Zhang2020} for convex minimization problems certifies the error by the primal--dual gap, which
equals the combined energy error with constant one, and uses balanced
adaptive greedy algorithms; the least-squares functional is discussed there
as a special case in section~6.5.  Chaudhry, Olson, and
Sentz~\cite{ChaudhryOlsonSentz2021} developed the least-squares finite
element RB method: the snapshots are least-squares finite element
solutions, and the reduced solution minimizes the least-squares functional
over their span.  They derived a rigorous bound of the error with respect
to the exact solution in the norm of $\cX$.  The bound uses a lower bound of
the parameter-dependent coercivity constant of the least-squares bilinear
form, computed by the successive constraint method, and the first-order
formulation is important for it.  In both~\cite{Yano2016}
and~\cite{ChaudhryOlsonSentz2021}, the snapshots are least-squares finite
element solutions, and a bound of the error in a norm that does not depend
on $\mu$ needs a lower bound of a stability constant.  Saddle-point problems
such as the Stokes equations are not considered there.

In section~3.3.2 of~\cite{QuarteroniManzoniNegri2016}, the term
``least-squares RB'' means the least-squares Petrov--Galerkin projection of
the truth problem.  Here, ``least-squares'' refers to the
least-squares functional of the first-order system of the partial
differential equation.  Ali, Steih, and Urban~\cite{AliSteihUrban2017}
compute the snapshots by adaptive wavelet methods and estimate the error of
the RB approximation with respect to the exact solution.  For the normal
equations, they show by a best approximation property that the reduced
residual at a selected parameter is not larger than the residual of its
snapshot.  In our setting, this is the second inequality
of~\eqref{eq:monotone}, where the residual is an $L^2$ residual and equals
the least-squares error of the whole pair.

\subsection{No inf--sup condition for the online reduced problem}\label{sec:noinfsup}
For indefinite problems, the Galerkin method is stable only
under an inf--sup condition.  For second-order elliptic equations, a
conforming finite element space satisfies it on sufficiently fine meshes, while a reduced space in general does not.  The
certification uses a lower bound of the inf--sup constant of the truth
operator, which can be computed by the successive constraint
method~\cite{HuynhEtAl2007,RozzaHuynhManzoni2013}.  For saddle-point
problems, a Galerkin RB method needs Brezzi's condition on the reduced
spaces, which can be restored by supremizers~\cite{RozzaVeroy2007}.  For
transport problems, the double greedy algorithm
of~\cite{DahmenPleskenWelper2014} builds a stabilizing test space together
with the trial space and controls the $L^2$ error.  In our method,
by~\ref{it:P1}, the online reduced problem is stable on any reduced space,
without any condition on the mesh size, and by~\ref{it:P2}, its certificate
needs no stability constant.

\subsection{Contributions}
Our method extends the least-squares RB methods of~\cite{Yano2016}
and~\cite{ChaudhryOlsonSentz2021} and the primal--dual RB method
of~\cite{Zhang2020} as follows.

\begin{enumerate}[label=(\roman*)]
\item By~\ref{it:P1}, the online reduced problem is the least-squares
projection onto the reduced space.  Since the lower bound in~\eqref{eq:normequiv} below holds
on the whole space, the online reduced problem is stable on any reduced space and for every
parameter, without any condition on the mesh size.  This holds for
indefinite second-order elliptic equations and for
saddle-point problems such as the Stokes equations: no discrete or reduced
inf--sup condition has to be verified, and no supremizers are needed.
\item By~\ref{it:P2}, the square root of the least-squares functional is the
norm of the residual of the continuous problem.  As the primal--dual gap in~\cite{Zhang2020}, it is the
certificate online and in the greedy algorithms, and no stability constant
is needed.  We do not convert it to a $\mu$-independent norm.
Such a conversion uses one number, the smallest singular value $\alpha(\mu)$
of $\cA_\mu:\cX\to\cY$, while $\tnorm{\cdot}_\mu$ uses the whole
operator; the bound
$\norm{U(\mu)-W}_{\cX}\le\eta(W;\mu)/\alpha(\mu)$ can overestimate the
error by any factor up to $\beta(\mu)/\alpha(\mu)$, with $\beta(\mu)$ the
largest singular value, even if $\alpha(\mu)$ is computed
exactly~\cite[Sect.~1.1]{Zhang2020}.  For a quantity of interest, the error
is a number, so there is no norm to choose, and no stability constant has
to be computed for each $\mu$;
see the companion paper~\cite{ZhangQoI2026}.
\item Since~\ref{it:P2} holds for every conforming pair, the snapshots can be
computed by any solver whose output can be made conforming by a
reconstruction.  The flux is recovered by least-squares, and the snapshot is
certified by its residual.  For the Stokes equations,
only the velocity of the black box is used.
\end{enumerate}

Section~\ref{sec:abstract} presents the method in an abstract setting, and
section~\ref{sec:greedy} the greedy algorithms.  Sections~\ref{sec:model},
\ref{sec:oseen}, and~\ref{sec:transport} treat second-order elliptic
equations, the Stokes equations, and linear transport, and
section~\ref{sec:numerics} reports the numerical experiments.

\section{The least-squares reduced basis method}
\label{sec:abstract}

\subsection{Abstract setting}

Let $\cX$ and $\cY$ be real Hilbert spaces, and let $\cD\subset\R^p$ be a
compact parameter domain.  For each $\mu\in\cD$, let $\cA_\mu:\cX\to\cY$ be
a bounded linear operator and $F_\mu\in\cY$.  We assume that
$\cA_\mu U(\mu)=F_\mu$ has a unique solution $U(\mu)\in\cX$.  An essential
boundary condition is either built into $\cX$ or imposed weakly as a
boundary residual, measured in an $L^2$ norm that is a component of $\cY$.
Nonhomogeneous essential data built into $\cX$ are first lifted.  The
reduced set is then a lifting plus a subspace of $\cX$, and the lifting at
the parameter of a snapshot is subtracted from the snapshot.  In all cases,
the functional is $\norm{\cA_\mu W-F_\mu}_{\cY}^2$, and~\eqref{eq:exactcert}
holds.

\begin{assum}[Uniform stability]\label{ass:uniform}
There exist constants $0<\aLS\le\bLS<\infty$, independent of $\mu\in\cD$,
such that
\begin{equation}\label{eq:normequiv}
  \aLS\norm{W}_{\cX}\ \le\ \norm{\cA_\mu W}_{\cY}\ \le\ \bLS\norm{W}_{\cX}
  \qquad\text{for all } W\in\cX,\ \mu\in\cD .
\end{equation}
\end{assum}

These constants are only used to compare the least-squares error with the
norm of $\cX$, for example, $\norm{U(\mu)-W}_{\cX}\le\aLS^{-1}\eta(W;\mu)$.
The certificate does not use them.

\subsection{The online reduced problem}

For a finite-dimensional subspace $X_N\subset\cX$, shifted by the lifting of
the essential data, the online reduced problem
gives the reduced solution and its certificate
\begin{equation}\label{eq:LSRB}
  U_N(\mu):=\arg\min_{W\in X_N}\norm{\cA_\mu W-F_\mu}_{\cY}^2
  \quad\text{and}\quad
  \eta_N(\mu):=\eta(U_N(\mu);\mu).
\end{equation}

\begin{prop}[Best approximation and exact certificate]\label{prop:best}
Under Assumption~\ref{ass:uniform}, $U_N(\mu)$ exists, is unique, and is the
orthogonal projection of $U(\mu)$ onto $X_N$ with respect to the inner
product $(\cA_\mu\cdot,\cA_\mu\cdot)_{\cY}$.  We have
\begin{equation}\label{eq:bestapprox}
\begin{gathered}
  \eta_N(\mu)=\tnorm{U(\mu)-U_N(\mu)}_\mu
  =\inf_{W\in X_N}\tnorm{U(\mu)-W}_\mu ,\\
  \norm{U(\mu)-U_N(\mu)}_{\cX}\le\frac{\bLS}{\aLS}\inf_{W\in X_N}\norm{U(\mu)-W}_{\cX}.
\end{gathered}
\end{equation}
If $X_N\subset X_M$, then $\eta_M(\mu)\le\eta_N(\mu)$ for all $\mu\in\cD$.
\end{prop}

\begin{proof}
By~\eqref{eq:exactcert}, the functional in~\eqref{eq:LSRB} is
$\tnorm{U(\mu)-W}_\mu^2$, and $\tnorm{\cdot}_\mu$ is a norm on $\cX$
by~\eqref{eq:normequiv}.  This gives the projection property and the first
relation in~\eqref{eq:bestapprox}.  The second one follows from
$\aLS\norm{U-U_N}_{\cX}\le\tnorm{U-U_N}_\mu=\inf_{W}\tnorm{U-W}_\mu
\le\bLS\inf_W\norm{U-W}_{\cX}$.  The monotonicity holds because the same
functional is minimized over a larger set.
\end{proof}

\subsection{Snapshots}

Let $\cX=\bSigma\times V$ be the product of a flux space and a potential
space, and let $\bSigma_h\subset\bSigma$ and $V_h\subset V$ be finite
element spaces on a mesh $\cT_h$ of our choice.  Denote
$\cJ_\mu(W):=\eta(W;\mu)^2$.  For a selected parameter
$\mu$, the user's solver returns an approximation $u_h(\mu)$ of the
potential, and a reconstruction $\cR_h$ maps it to
$\widetilde u_h(\mu)\in V_h$.  If the solver is conforming on $\cT_h$,
$\cR_h$ is the identity.  The flux is recovered by
\begin{equation}\label{eq:recovery}
  \bsigma_h(\mu):=\arg\min_{\btau\in\bSigma_h}
  \cJ_\mu\bigl(\btau,\widetilde u_h(\mu)\bigr).
\end{equation}
This is the flux block of the least-squares finite element method and is
solved once for each snapshot.  If the solver provides a conforming flux,
this flux can be used instead.

\begin{prop}[Best recovery]\label{prop:recbest}
Let $\bsigma_h(\mu)$ be the minimizer in~\eqref{eq:recovery}, and let
$\widetilde U_h(\mu):=(\bsigma_h(\mu),\widetilde u_h(\mu))$.  Then
\[
  \eta\bigl(\widetilde U_h(\mu);\mu\bigr)
  =\min_{\btau\in\bSigma_h}\tnorm{U(\mu)-(\btau,\widetilde u_h(\mu))}_\mu .
\]
\end{prop}

\begin{proof}
The proof is the same as that of Proposition~\ref{prop:best}, with the
affine set $\bSigma_h\times\{\widetilde u_h(\mu)\}$ in place of $X_N$.
\end{proof}

\begin{rem}\label{rem:primaryfloor}
For every $\btau$,
\[
  \tnorm{U(\mu)-(\btau,\widetilde u_h)}_\mu
  \ge\aLS\norm{U(\mu)-(\btau,\widetilde u_h)}_{\cX}
  \ge\aLS\norm{u(\mu)-\widetilde u_h(\mu)}_V .
\]
Hence the flux recovery cannot make up for an inaccurate reconstructed
potential, and a poor solver or a poor reconstruction leads to a large
snapshot certificate.
\end{rem}

The snapshot is certified by
\begin{equation}
  \eta_h^{\snap}(\mu):=\eta\bigl(\widetilde U_h(\mu);\mu\bigr)
  =\tnorm{U(\mu)-\widetilde U_h(\mu)}_\mu .
\end{equation}
The certificate does not use any property of $\cR_h$.  The quality of
$\cR_h$ and of the solver only affects the value of $\eta^\snap_h$.  The
snapshots at different $\mu_i$ may come from different solvers and
different meshes.

\subsection{Reduced spaces and offline--online decomposition}
\label{sec:offline}

For the selected parameters $\mu_1,\dots,\mu_N$, let
$\widetilde U_i=(\widetilde\bsigma_i,\widetilde u_i):=\widetilde U_h(\mu_i)$
and define
\begin{equation}
  \bSigma_N:=\Span\{\widetilde\bsigma_i\}_{i=1}^N,\qquad
  V_N:=\Span\{\widetilde u_i\}_{i=1}^N,\qquad
  X_N:=\bSigma_N\times V_N .
\end{equation}
The online reduced problem uses the product space $X_N$, whose dimension is at most
$2N$, instead of the span of the snapshot pairs.  Thus the fluxes and the
potentials have independent coefficients, as in the decoupled formulation
of~\cite[Sect.~2.4]{Yano2016}.  For
conditioning, the two families are orthonormalized separately by the
Gram--Schmidt process with reorthogonalization, in fixed inner products of
$\bSigma$ and $V$.  We denote the orthonormal bases by
$\{\btau_j\}_{j=1}^{N_\sigma}$ and $\{v_j\}_{j=1}^{N_u}$ and their union by
$\{\varphi_j\}$.  For $W=\sum_jx^\sigma_j\btau_j+\sum_jx^u_jv_j$ and
$x=(x^\sigma,x^u)$, we have
\begin{equation}\label{eq:quadform}
\begin{gathered}
  \cJ_\mu(W)=x^TM(\mu)x-2x^Tg(\mu)+\cJ_\mu(0),\\
  M_{ij}(\mu)=(\cA_\mu\varphi_j,\cA_\mu\varphi_i)_{\cY},\qquad
  g_i(\mu)=(F_\mu,\cA_\mu\varphi_i)_{\cY},
\end{gathered}
\end{equation}
and $U_N(\mu)$ is given by $M(\mu)x=g(\mu)$.  If the forms
$(\cA_\mu W,\cA_\mu W')_{\cY}$ and $(F_\mu,\cA_\mu W)_{\cY}$ and the number
$\cJ_\mu(0)=\norm{F_\mu}_{\cY}^2$ are parameter-separable, then the parts of $M(\mu)$,
$g(\mu)$, and $\cJ_\mu(0)$ that do not depend on $\mu$ are assembled offline, and the online cost does not depend on the dimension of
the truth spaces.  The computed certificate is exact if the
integrals in~\eqref{eq:quadform} are computed exactly, which is the case in
all our experiments.  If approximate data $\widetilde F_\mu$ are used, the
certificate changes by at most $\norm{F_\mu-\widetilde F_\mu}_{\cY}$.

\subsection{Snapshots on different meshes}
\label{sec:assembly}

The offline matrices need integrals of products of snapshots computed on
different meshes.  We only consider the case of an adaptive
greedy algorithm, where the snapshot meshes are nested,
$\cT_1\subset\cT_2\subset\dots\subset\cT_N$, and each mesh is a conforming
refinement of the previous one.  The finite element spaces of fixed degree
used here (Lagrange, Raviart--Thomas, and discontinuous) are nested under
refinement, so a snapshot on $\cT_i$ can be represented exactly on $\cT_j$
for $j\ge i$.  The entries between the new snapshot on $\cT_N$ and an old
snapshot are therefore computed on $\cT_N$.  For coefficients that are
piecewise constant on the initial
mesh, these integrals are exact.  Each new snapshot adds one row and one
column for each block and each separable term.
\section{Greedy algorithms with residual certificates}\label{sec:greedy}

The certificate is the norm of the residual of the continuous problem.  On a fixed truth
space, it cannot be smaller than the least-squares residual in this space, so
the classical greedy algorithm does not stop if the RB tolerance is below
this residual.  The algorithms in section~7
of~\cite{Zhang2020} balance the truth and the RB tolerances to solve this
problem.  In this section, we present their non-intrusive versions.  Let
$\Xi_\train\subset\cD$ be a finite training set, let $\varepsilon_h$ be the
truth tolerance for the snapshot certificates and $\varepsilon_\rb$ the RB
tolerance for the reduced certificates, let $\rrf>1$ be their ratio, and let
$N_{\max}$ be a maximal dimension.

\subsection{Non-intrusive level 0: least-squares truth}

Suppose that the truth solver is the least-squares method on a finite
element space $X_h\subset\cX$.  Then the snapshot at $\mu_i$ is the
projection $U_h(\mu_i):=\arg\min_{W\in X_h}\eta(W;\mu_i)$, with the
certificate $\eta_h(\mu_i):=\eta(U_h(\mu_i);\mu_i)$.  On a fixed mesh, any
reduced space spanned by such snapshots satisfies $X_N\subset X_h$.  By
Proposition~\ref{prop:best}, we obtain the two results on which the greedy
algorithms of~\cite{Zhang2020} are based: $\eta_N(\mu)\ge\eta_h(\mu)$ for
all $\mu$, and $\eta_N(\mu_i)=\eta_h(\mu_i)$ at the selected parameters,
since $U_h(\mu_i)\in X_N$; see also Theorem~5.7 of~\cite{Zhang2020}.  In
fact, $U_h(\mu)$ is the orthogonal projection of $U(\mu)$ onto $X_h$ with
respect to the inner product in~\ref{it:P1}, and $X_N\subset X_h$, so
\begin{equation}\label{eq:split}
  \eta_N(\mu)^2=\eta_h(\mu)^2+\tnorm{U_h(\mu)-U_N(\mu)}_\mu^2
  \qquad\text{for all }\mu .
\end{equation}
The second term is the error with respect to the truth solution, which is
the quantity that standard RB error bounds certify for a fixed truth
discretization.  For frozen snapshots on nested meshes (see below),
\eqref{eq:split} holds with $X_h$ the space on the finest mesh.  At
non-intrusive levels~1--3, the snapshot is not a projection.  The reconstructed pairs are
in $X_h$, so $\eta_N\ge\eta_h$ and~\eqref{eq:split} still hold, where
$U_h(\mu)$ is the least-squares solution in $X_h$ (it is not computed).  The
equality at the selected parameters becomes the inequality
$\eta_N(\mu_i)\le\eta^{\snap}_h(\mu_i)$.

\subsection{Fixed resolution}

Algorithm~\ref{alg:fixed} is Algorithm~1 of~\cite{Zhang2020} with the
snapshot certificate in place of the finite element error estimator.
By~\eqref{eq:split}, $\max_{\Xi_\train}\eta_N\ge\max_{\Xi_\train}\eta_h$,
where the right-hand side is the largest certificate of the least-squares
solution on the mesh, and $\eta_h(\mu_i)\le\eta_N(\mu_i)\le
\eta^\snap_h(\mu_i)$ at a selected parameter.  The RB tolerance follows the
largest snapshot certificate $\varepsilon^N_h$.
\par\medskip\noindent\begin{minipage}{\textwidth}
\begin{algo}[Fixed resolution, adaptive RB tolerance]\label{alg:fixed}
Input: $\Xi_\train$, a fixed mesh (or grid), $\rrf>1$, $\varepsilon_\rb^0$,
$N_{\max}$, $\mu_1$; $\varepsilon_h^0=0$, $N=1$.
\begin{enumerate}[label=\arabic*:]
\item Call the user's solver at $\mu_N$; reconstruct, recover the flux, and
compute $\eta^\snap_h(\mu_N)$.
\item Let $\varepsilon_h^N=\max\{\eta^\snap_h(\mu_N),\varepsilon_h^{N-1}\}$
and $\varepsilon_\rb^N=\max\{\rrf\varepsilon_h^N,\varepsilon_\rb^{N-1}\}$.
\item Update $\bSigma_N$, $V_N$, and the offline matrices.
\item Choose $\mu_{N+1}=\arg\max_{\mu\in\Xi_\train}\eta_N(\mu)$.
\item If $\eta_N(\mu_{N+1})>\varepsilon_\rb^N$ and $N<N_{\max}$, set
$N:=N+1$ and go to~1; otherwise, stop.
\end{enumerate}
\end{algo}
\end{minipage}\par\medskip

\subsection{Adaptive resolution with frozen snapshots}

Algorithm~\ref{alg:adaptive} is Algorithm~2 of~\cite{Zhang2020} with frozen
snapshots: as in~\cite[Sect.~4.2]{Yano2016}, the earlier snapshots are kept
and represented on the refined mesh, not recomputed.  It requires a solver that accepts a mesh;
a solver on its own grid is used with Algorithm~\ref{alg:fixed}.  Since the certificate is exact for every
conforming pair, it remains exact when the snapshots come from different meshes.
The reduced spaces are nested, so, as with the residual bound
of~\cite[Sect.~2.4]{Yano2016}, the certificate is monotone in $N$, and each
frozen snapshot already meets the tolerance $\varepsilon_h$ at its own
parameter.
\par\medskip\noindent\begin{minipage}{\textwidth}
\begin{algo}[Adaptive resolution, fixed tolerances, frozen snapshots]
\label{alg:adaptive}
Input: $\Xi_\train$, $\varepsilon_h$, $\rrf>1$,
$\varepsilon_\rb=\rrf\varepsilon_h$, $N_{\max}$, an initial mesh $\cT_0$,
$\mu_1$; $N=1$.
\begin{enumerate}[label=\arabic*:]
\item Starting from $\cT_{N-1}$, call the user's solver at $\mu_N$,
reconstruct, and recover the flux.  While
$\eta^\snap_h(\mu_N)>\varepsilon_h$, mark the elements by the local
contributions of $\cJ_{\mu_N}$, refine, and repeat.  The final mesh is
$\cT_N$.
\item Freeze $\widetilde U_h(\mu_N)$ on $\cT_N$.  The earlier snapshots are
not recomputed.
\item Update $\bSigma_N$, $V_N$, and the offline matrices
(section~\ref{sec:assembly}).
\item Choose $\mu_{N+1}=\arg\max_{\mu\in\Xi_\train}\eta_N(\mu)$.
\item If $\eta_N(\mu_{N+1})>\varepsilon_\rb$ and $N<N_{\max}$, set $N:=N+1$
and go to~1; otherwise, stop.
\end{enumerate}
\end{algo}
\end{minipage}\par\medskip
Step~1 assumes that $\eta^\snap_h(\mu_N)$ tends to zero under this
refinement; at non-intrusive level~1, this follows from the convergence result
of~\cite{LiZhang2025}.

\begin{lem}[Monotonicity and no re-selection]\label{lem:monotone}
Suppose that the snapshots are frozen.  Then
$X_{N-1}\subset X_N$ and
\begin{equation}\label{eq:monotone}
  \eta_N(\mu)\le\eta_{N-1}(\mu)\quad\text{for all }\mu,
  \qquad
  \eta_N(\mu_i)\le\eta^\snap_h(\mu_i)\quad(1\le i\le N).
\end{equation}
Therefore, $\eta_N(\mu_i)\le\varepsilon^N_h<\varepsilon^N_\rb$ in
Algorithm~\ref{alg:fixed} and $\eta_N(\mu_i)\le\varepsilon_h<\varepsilon_\rb$
in Algorithm~\ref{alg:adaptive}, and no parameter is selected twice.
\end{lem}

\begin{proof}
Both inequalities hold because the same functional is minimized over a
larger set (Proposition~\ref{prop:best}): $X_N\supset X_{N-1}$ and
$\widetilde U_i\in X_N$.
\end{proof}

\begin{rem}\label{rem:relcert}
If the zero pair is admissible (homogeneous essential data or weakly imposed
boundary data), then $\cJ_\mu(0)^{1/2}=\tnorm{U(\mu)}_\mu$
by~\eqref{eq:exactcert}.  Hence the relative error
$\tnorm{U(\mu)-U_N(\mu)}_\mu/\tnorm{U(\mu)}_\mu$ is
$r_N(\mu):=\eta_N(\mu)/\cJ_\mu(0)^{1/2}$.  The number $\cJ_\mu(0)$ is
the last term in~\eqref{eq:quadform}.  It is available online at no extra
cost, and it depends neither on $N$ nor on the mesh if the data terms are
integrated exactly.  Dividing~\eqref{eq:monotone} by $\cJ_\mu(0)^{1/2}$, we
obtain $r_N(\mu)\le r_{N-1}(\mu)$ and
$r_N(\mu_i)\le\eta^\snap_h(\mu_i)/\cJ_{\mu_i}(0)^{1/2}$.  Thus
Lemma~\ref{lem:monotone} and Proposition~\ref{prop:sat} also hold for
$r_N$, and the algorithms can be run with $r_N$ in place of $\eta_N$, with
the snapshot test $\eta^\snap_h(\mu_N)\le\varepsilon_h\cJ_{\mu_N}(0)^{1/2}$
and relative tolerances.  This is useful when $\tnorm{U(\mu)}_\mu$ varies by
orders of magnitude over $\cD$.  For~\eqref{eq:cdtest}, $\cJ_\mu(0)=1$ for
all $\mu$.
\end{rem}

\subsection{Saturation, enriched training sets, and the safety check}
In step~4, the certificate is evaluated on the whole training set.  The
saturation trick~\cite[Def.~3.3, Alg.~2]{HesthavenStammZhang2014} skips a
parameter if its stored estimate, multiplied by a saturation constant, is
smaller than the current maximum.  For a general error estimator, this is an
assumption, and the user chooses the saturation constant.  The first inequality in~\eqref{eq:monotone} is the saturation
property with constant one, as on a fixed mesh
in~\cite[Sect.~7.2]{Zhang2020}.

\begin{prop}[Saturation and removal]\label{prop:sat}
Suppose that the snapshots are frozen.
\begin{enumerate}[label=(\alph*)]
\item If $n<N$ and $\eta_n(\mu)\le m$, then $\eta_N(\mu)\le m$.  Hence, with
saturation constant one, the maximum of $\eta_N$ over the training set is
attained at a parameter whose certificate is evaluated, and the saturation
trick selects the same parameter as the full search.
\item In Algorithms~\ref{alg:fixed} and~\ref{alg:adaptive} (with
$\varepsilon_\rb^N:=\varepsilon_\rb$ in the latter), if
$\eta_n(\mu)\le\varepsilon_\rb^n$, then $\eta_N(\mu)\le\varepsilon_\rb^N$
for all $N\ge n$.  Hence a parameter removed from the training set for this
reason would never be selected.
\end{enumerate}
\end{prop}

\begin{proof}
We have $\eta_N(\mu)\le\eta_n(\mu)$ by~\eqref{eq:monotone}, and
$\varepsilon_\rb^N$ is nondecreasing in $N$.
\end{proof}

\begin{rem}
The adaptively enriched greedy algorithm
of~\cite[Sect.~3.3, Alg.~3]{HesthavenStammZhang2014} selects parameters
from a small random training set and replaces the parameters below the
tolerance by new ones.  By Proposition~\ref{prop:sat}, it can be used with
saturation constant one.  Its final safety check uses $n$ new random
parameters.  For $\varepsilon>0$, if the set $\{\mu:\eta_N(\mu)>\varepsilon\}$ has probability
at least $p$, a reduced space passes this check with probability at most
$(1-p)^n\le e^{-pn}$.  By~\eqref{eq:exactcert}, this is the set of parameters at
which the error with respect to the exact solution is larger than $\varepsilon$.
\end{rem}

\Needspace{6\baselineskip}
\section{Second-order elliptic equations}
\label{sec:model}

\subsection{The model problem}

Let $\Omega\subset\R^d$ be a bounded Lipschitz polyhedral domain with
$\partial\Omega=\overline{\Gamma_D}\cup\overline{\Gamma_N}$ and $|\Gamma_D|>0$.  For $\mu\in\cD$,
consider
\begin{equation}
\begin{gathered}
  -\divv(A_\mu\grad u)+\bb_\mu\cdot\grad u+c_\mu u=f_{1,\mu}+\divv\bff_{2,\mu}
  \ \text{ in }\Omega,\\
  u=\phi_\mu\ \text{ on }\Gamma_D,\qquad
  -(A_\mu\grad u+\bff_{2,\mu})\cdot\bn=\psi_\mu\ \text{ on }\Gamma_N.
\end{gathered}
\end{equation}
We assume that $A_\mu\in L^\infty(\Omega)^{d\times d}$ is symmetric and
uniformly positive definite, that $\bb_\mu\in L^\infty(\Omega)^d$ and
$c_\mu\in L^\infty(\Omega)$, with bounds independent of $\mu$, and that the
bilinear form
$a_\mu(w,v):=(A_\mu\grad w,\grad v)+(\bb_\mu\cdot\grad w+c_\mu w,v)$
satisfies the inf--sup condition on $H^1_D(\Omega)$ with a constant
$\beta_0>0$ independent of $\mu\in\cD$.  With
$\bsigma=-(A_\mu\grad u+\bff_{2,\mu})$, the first-order system is
\begin{equation}
  \divv\bsigma+\bb_\mu\cdot\grad u+c_\mu u=f_{1,\mu}\quad\text{and}\quad
  \bsigma+A_\mu\grad u=-\bff_{2,\mu}\quad\text{in }\Omega,
\end{equation}
with the two boundary conditions, and the volume functional is
\begin{equation}\label{eq:Jelliptic}
  \LS_\mu(\btau,v)
  :=\norm{A_\mu^{-1/2}(\btau+A_\mu\grad v+\bff_{2,\mu})}_{0}^{2}
   +\norm{\divv\btau+\bb_\mu\cdot\grad v+c_\mu v-f_{1,\mu}}_{0}^{2},
\end{equation}
where $\norm{\cdot}_0$ and $\norm{\cdot}_1$ denote the norms of $L^2(\Omega)$ and
$H^1(\Omega)$.  This is $\cJ_\mu$ of section~\ref{sec:abstract} with
$\cY=L^2(\Omega)^d\times L^2(\Omega)$,
$\cA_\mu(\btau,v)=\bigl(A_\mu^{-1/2}(\btau+A_\mu\grad v),\
\divv\btau+\bb_\mu\cdot\grad v+c_\mu v\bigr)$, and
$F_\mu=\bigl(-A_\mu^{-1/2}\bff_{2,\mu},\ f_{1,\mu}\bigr)$.
Let $\cX:=H_N(\divvr;\Omega)\times H^1_D(\Omega)$, where the subscripts denote zero
normal trace on $\Gamma_N$ and zero trace on $\Gamma_D$.  Then
Assumption~\ref{ass:uniform} holds on $\cX$ in the norm
$(\norm{\btau}_{\Hdiv}^2+\norm{v}_1^2)^{1/2}$, with constants depending on
$\beta_0$ and the bounds of the
coefficients~\cite{CaiLazarovManteuffelMcCormick1994,Zhang2023}.

The right-hand side $f_{1,\mu}+\divv\bff_{2,\mu}$ is a general element of
$H^{-1}(\Omega)$, since every element of $H^{-1}(\Omega)$ has this form with
$f_1\in L^2(\Omega)$ and $\bff_2\in L^2(\Omega)^d$.  The vector part is included in
the flux $\bsigma=-(A_\mu\grad u+\bff_{2,\mu})$, which is an exact change of
variables.  Then $\divv\bsigma\in L^2(\Omega)$, and~\eqref{eq:Jelliptic} only
contains $L^2$ norms; no discrete $H^{-1}$ inner product as
in~\cite{BLP1997} is used.  This treatment of $H^{-1}$ data is that
of~\cite{Zhang2023}; see also~\cite{WuZhang2026}.  For the RB method,
parameter-dependent $H^{-1}$ data, for example line sources, enter the
certificate, the flux recovery, and the affine decomposition in the same
way as $L^2$ data.

\subsection{Boundary conditions}
Suppose that $\phi_\mu$ and $\psi_\mu$ are satisfied exactly by every
reduced pair.  This is the case for homogeneous data and for data that are
discrete on the coarsest mesh and enter through a parameter-separable
lifting as in~\cite{Zhang2020}.  Then $\cX$ is the space with homogeneous
essential conditions after the lifting, and~\eqref{eq:exactcert} holds with
$\tnorm{\cdot}_\mu=\LS_\mu^{0}(\cdot)^{1/2}$, where $\LS_\mu^{0}$ is the
functional with zero data.  If the reduced pairs satisfy discrete data
$\phi_{h,\mu}$ and $\psi_{h,\mu}$ exactly, the certificate is exact with
respect to the solution of the continuous problem with these data.
Independently of any discretization, the
distance between this solution and $U(\mu)$ in the norm of
$\Hdiv\times H^1(\Omega)$ is equivalent to
$\norm{\phi_\mu-\phi_{h,\mu}}_{1/2,\Gamma_D}+\norm{\psi_\mu-\psi_{h,\mu}}_{\Hmoo(\Gamma_N)}$,
where $\Hmoo(\Gamma_N)$ is the dual of $H^{1/2}_{00}(\Gamma_N)$.  This term depends
on $\mu$ and on the boundary mesh but not
on $N$.

\subsection{Snapshots and offline--online decomposition}
\label{sec:reconstruction}

Let $V_h$ be a continuous piecewise polynomial space of degree $k$ and
$\bSigma_h$ the Raviart--Thomas space of index $k-1$ on $\cT_h$.  The reconstruction $\cR_h$ maps the
output of the user's solver to $V_h$.
\begin{itemize}
\item \emph{Non-intrusive level 1.}  For a conforming Galerkin solution in $V_h$ on
$\cT_h$, $\cR_h$ is the identity.
\item \emph{Non-intrusive level 2.}  For a discontinuous Galerkin or nonconforming
solution on $\cT_h$, $\cR_h$ averages the values at each Lagrange node over
the elements sharing the node, which gives an element of
$V_h$~\cite{KarakashianPascal2003}.  If the coefficient jumps across
interfaces, we use the modified Oswald interpolation of~\cite{CaiYeZhang2011}:
at a vertex on an interface, the values are averaged only over the elements
of the subdomain with the largest coefficient.  For a cell-centered finite
volume solution, the nodal values can be taken as averages of the values of
the adjacent cells.
\item \emph{Non-intrusive level 3.}  For a finite difference solution on a grid, $\cR_h$
uses the grid values as the nodal values of a continuous piecewise bilinear
function on the cells of the grid, or of a continuous piecewise linear
function on a triangulation of the grid.
\end{itemize}

\begin{rem}
Let $u_h$ be a discontinuous Galerkin or Crouzeix--Raviart approximation of
$u(\mu)$ with homogeneous Dirichlet data, let $\widetilde u_h=\cR_hu_h$ be
its nodal average, and let $\norm{\cdot}_{\rm DG}$ be the broken $H^1$ norm
with the jump term $\sum_Fh_F^{-1}\norm{[\![v]\!]}_{0,F}^2$, where boundary
edges are included.  The jumps of $u_h$ are those of $u_h-u$.  By
Proposition~\ref{prop:recbest}, the upper bound in~\eqref{eq:normequiv},
the averaging estimate of~\cite{KarakashianPascal2003} (of~\cite{CaiYeZhang2011}
for the modified Oswald interpolation), and the Poincar\'e inequality, we
have
\[
  \eta^\snap_h(\mu)\le\bLS\Bigl(\inf_{\btau\in\bSigma_h}
  \norm{\bsigma(\mu)-\btau}_{\Hdiv}+C'\norm{u(\mu)-u_h}_{\rm DG}\Bigr),
\]
where $C'$ depends only on the shape regularity, the polynomial degree, and
the Poincar\'e constant.  By Remark~\ref{rem:primaryfloor}, we also have
$\eta^\snap_h(\mu)\ge\aLS\norm{u(\mu)-\widetilde u_h}_1$.  Hence a black box
that is accurate in its own norm gives a snapshot with a small certificate.
Certifying $u_h$ itself would require the jump term, which depends on the
solver; the snapshot certificate does not.
\end{rem}
Every snapshot must satisfy the discrete data assumed by the online
lifting.  The reconstruction $\cR_h$ sets the boundary nodal values of the
potential; this replaces the boundary values of a black box that imposes
the Dirichlet condition weakly.  The flux recovery~\eqref{eq:recovery} is
over the Raviart--Thomas functions satisfying the discrete Neumann data, and
its normal equations are
\begin{equation}
\begin{aligned}
  (A_\mu^{-1}\bsigma_h,\bchi)+(\divv\bsigma_h,\divv\bchi)
  &=-(\grad\widetilde u_h+A_\mu^{-1}\bff_{2,\mu},\bchi)\\
  &\quad+(f_{1,\mu}-\bb_\mu\cdot\grad\widetilde u_h-c_\mu\widetilde u_h,\divv\bchi)
\end{aligned}
\end{equation}
for all such $\bchi$ with zero Neumann data~\cite{CaiZhang2010,LiZhang2025}.
If the solver accepts a mesh, the local contributions of $\LS_\mu$ at
$\widetilde U_h(\mu)$ are used as refinement indicators.  At non-intrusive level~1, the
adaptive algorithm driven by these indicators
converges~\cite{LiZhang2025}.

\section{The Stokes equations without reduced inf--sup conditions}
\label{sec:oseen}

\subsection{Stress--velocity formulation}
Let $\Omega\subset\R^d$ be a bounded Lipschitz polyhedral domain, let
$\nu_\mu$ be piecewise constant on a fixed partition of $\Omega$ with
$0<\nu_{\min}\le\nu_\mu\le\nu_{\max}$, and consider
\begin{equation}\label{eq:oseen}
  -\divv\bigl(\nu_\mu\beps(\bu)\bigr)+\grad p=\bff_\mu,
  \quad \divv\bu=0\ \text{ in }\Omega,\qquad
  \bu=\bzero\ \text{ on }\partial\Omega,\quad p\in L^2_0(\Omega),
\end{equation}
where $\beps(\bu)=\frac12(\grad\bu+(\grad\bu)^t)$ and
$\bff_\mu\in L^2(\Omega)^d$.  For a discontinuous viscosity, the symmetric form
$\divv(\nu_\mu\beps(\bu))$ gives the physical interface condition, the
continuity of the traction, while the pseudostress $\nu_\mu\grad\bu-p\bI$
does not~\cite[Rems.~1 and~3]{LiangZhang2026}.  Therefore, we use the
stress.  Let $\dev\btau:=\btau-\frac1d(\tr\btau)\bI$ and
$\bsigma:=\nu_\mu\beps(\bu)-p\bI$.  Then $\divv\bsigma=-\bff_\mu$ and
$\tr\bsigma=-dp$, and since $\divv\bu=0$, we have
$\dev\bsigma=\nu_\mu\beps(\bu)$.  Hence
\begin{equation}\label{eq:oseenfos}
  \nu_\mu^{-1/2}\dev\bsigma-\nu_\mu^{1/2}\beps(\bu)=\bzero
  \quad\text{and}\quad \divv\bsigma+\bff_\mu=\bzero\quad\text{in }\Omega,
\end{equation}
with $\bu=\bzero$ on $\partial\Omega$, and $p=-\frac1d\tr\bsigma$.  The condition
$p\in L^2_0(\Omega)$ becomes $\tr\bsigma\in L^2_0(\Omega)$.  Let
$\bSigma_0:=\{\btau\in\Hdiv^d:\tr\btau\in L^2_0(\Omega)\}$, where $\Hdiv^d$ is
taken row-wise and no symmetry is imposed, and let
$\cX:=\bSigma_0\times H^1_0(\Omega)^d$.  For $(\btau,\bv)\in\cX$, define
\begin{equation}\label{eq:oseenLS}
  \LS_\mu(\btau,\bv):=
  \norm{\nu_\mu^{-1/2}\dev\btau-\nu_\mu^{1/2}\beps(\bv)}_0^2
  +\norm{\nu_\mu^{-1/2}\bigl(\divv\btau+\bff_\mu\bigr)}_0^2 ,
\end{equation}
and let $\LS^0_\mu$ be the same functional with zero data.  This is
$\cJ_\mu$ of section~\ref{sec:abstract} with
$\cY=L^2(\Omega)^{d\times d}\times L^2(\Omega)^d$,
$\cA_\mu(\btau,\bv)=\bigl(\nu_\mu^{-1/2}\dev\btau-\nu_\mu^{1/2}\beps(\bv),\
\nu_\mu^{-1/2}\divv\btau\bigr)$, and
$F_\mu=\bigl(\bzero,\ -\nu_\mu^{-1/2}\bff_\mu\bigr)$.  It is the
functional $J_1$ in~\cite[(7.1)]{LiangZhang2026}, which is defined there on
the subspace $\nu_\mu^{-1}\tr\btau\in L^2_0(\Omega)$.  This subspace depends on
$\mu$, while ours does not.  The map $\btau\mapsto\btau+c\bI$ between the two
subspaces changes neither $\dev\btau$ nor $\divv\btau$, so the norm
equivalence of~\cite{CaiLeeWang2004,LiangZhang2026} also holds on $\cX$.
Thus Assumption~\ref{ass:uniform} holds with
$\norm{(\btau,\bv)}_{\cX}^2:=\norm{\btau}_0^2+\norm{\divv\btau}_0^2+\norm{\grad\bv}_0^2$
and with constants depending on $d$, $\Omega$, $\nu_{\min}$, and $\nu_{\max}$.
The certificate does not use these constants.

\subsection{The certificate, the divergence, the symmetry, and the pressure}

\begin{prop}\label{prop:oseen}
Let $(\bsigma,\bu)\in\cX$ be the solution of~\eqref{eq:oseenfos}.  For
$(\btau,\bv)\in\cX$, let $\eta:=\LS_\mu(\btau,\bv)^{1/2}$,
$p_{\btau}:=-\frac1d\tr\btau$, and $\skw\bM:=\frac12(\bM-\bM^t)$.  Then
(a)~$\eta=\LS^0_\mu(\bsigma-\btau,\bu-\bv)^{1/2}$;
(b)~$d^{-1}\norm{\nu_\mu^{1/2}\divv\bv}_0^2
+\norm{\nu_\mu^{-1/2}\skw\btau}_0^2\le\eta^2$;
(c)~$\norm{p-p_{\btau}}_0\le d^{-1/2}\norm{\bsigma-\btau}_0$.
\end{prop}

\begin{proof}
(a) This is~\eqref{eq:exactcert}, since both residuals are affine and
vanish at $(\bsigma,\bu)$.  (b) In the Frobenius inner product, a matrix
$\bM\in\R^{d\times d}$ is the orthogonal sum of its symmetric deviatoric
part, $\skw\bM$, and $\frac1d(\tr\bM)\bI$, and the norm of the last part is
$d^{-1/2}|\tr\bM|$.  Let $\bR_1$ be the first residual
in~\eqref{eq:oseenLS}.  Since $\beps(\bv)$ and $\bI$ are symmetric,
$\skw\bR_1=\nu_\mu^{-1/2}\skw\btau$ and
$\tr\bR_1=-\nu_\mu^{1/2}\divv\bv$.  (c) We have
$p-p_{\btau}=-\frac1d\tr(\bsigma-\btau)$ and $|\tr\bM|\le\sqrt d\,|\bM|$.
\end{proof}

By~(b), the divergence of the velocity and the skew-symmetric part of the
stress are bounded by the certificate, for the black-box snapshots as well
as for the reduced solution; neither is assumed to vanish.  By~(c) and the norm
equivalence, the pressure error is also bounded by the certificate, with a
constant depending on $\nu_{\min}$ and $\nu_{\max}$.

\subsection{Non-intrusive snapshots}
The black box returns a velocity $\bu_h(\mu)$; a pressure, if computed, is
not used.  Let
$\widetilde\bu_h=\cR_h\bu_h\in V_h^d$, where $V_h$ consists of continuous
piecewise polynomials vanishing on $\partial\Omega$.  The reconstruction $\cR_h$ is
the identity for Taylor--Hood or stabilized equal-order $P_1$--$P_1$
velocities, and it drops the bubbles of a MINI velocity.  For nonconforming
or discontinuous Galerkin velocities, $\cR_h$ is the nodal average, where at
a vertex on an interface of $\nu_\mu$ only the elements of the subdomain
with the largest $\nu_\mu$ are used, as in section~\ref{sec:reconstruction}.
For the staggered (MAC) finite difference scheme on a Cartesian grid, $\cR_h$
takes at each grid vertex the mean of the two neighboring face values of
each velocity component (zero on $\partial\Omega$), and uses $P_1$ on the grid cut
into triangles.  Let $\bSigma_h^d$ be the row-wise Raviart--Thomas space of
index $k-1$ for velocities of degree $k$, for example, $RT_1$ for
Taylor--Hood.
The stress is recovered by
$\bsigma_h(\mu):=\arg\min_{\btau\in\bSigma_h^d\cap\bSigma_0}
\LS_\mu(\btau,\widetilde\bu_h(\mu))$, that is, for all
$\bchi\in\bSigma_h^d\cap\bSigma_0$,
\begin{equation}\label{eq:oseenrec}
  (\nu_\mu^{-1}\dev\bsigma_h,\dev\bchi)+(\nu_\mu^{-1}\divv\bsigma_h,\divv\bchi)
  =\bigl(\beps(\widetilde\bu_h),\dev\bchi\bigr)-(\nu_\mu^{-1}\bff_\mu,\divv\bchi).
\end{equation}
The bilinear form on the left is positive definite on $\bSigma_0$.  Indeed,
if it vanishes at $\bchi$, then $\dev\bchi=\bzero$ and $\divv\bchi=\bzero$,
so $\bchi=q\bI$ with $\grad q=\bzero$, and $q=0$ since
$\tr\bchi\in L^2_0(\Omega)$.  Both sides of~\eqref{eq:oseenrec} vanish at
$\bchi=\bI$, and the constant multiples of $\bI$ are in $\bSigma_h^d$, so
any solution of~\eqref{eq:oseenrec} on $\bSigma_h^d$ gives $\bsigma_h$ after
a multiple of $\bI$ is subtracted to make the mean trace zero.  The recovered
pressure is $p_h=-\frac1d\tr\bsigma_h$, and the snapshot certificate is
$\eta^\snap_h(\mu)=\LS_\mu(\bsigma_h,\widetilde\bu_h)^{1/2}$.  The recovery
uses only the velocity of the black box.  Thus the black box need not be inf--sup stable, and its velocity need not
be divergence free; its pressure space and its stabilization only affect the
value of $\eta^\snap_h$.

\subsection{The online reduced problem}
Let $V_N=\Span\{\widetilde\bu_h(\mu_i)\}\subset H^1_0(\Omega)^d$ and
$\bSigma_N=\Span\{\bsigma_h(\mu_i)\}\subset\bSigma_0$.  The online reduced problem
minimizes $\LS_\mu$ over $\bSigma_N\times V_N$; for its minimizer
$(\bsigma_N,\bu_N)$, the reduced pressure is
$p_N=-\frac1d\tr\bsigma_N$.  The reduced problem needs no pressure space and no supremizers
(section~\ref{sec:noinfsup}), and
Proposition~\ref{prop:oseen} applies to the reduced solution.  For $\nu_\mu$ piecewise constant on a fixed
partition, the offline--online decomposition of section~\ref{sec:offline}
can be used.  The frozen snapshots and the greedy algorithms of
section~\ref{sec:greedy} are used without change.

\Needspace{6\baselineskip}
\section{Linear transport with a moving inflow boundary}
\label{sec:transport}

\subsection{The flux formulation}
Let $\Omega\subset\R^d$ be a bounded Lipschitz polyhedral domain, let
$\bbeta_\mu\in\R^d$ be a constant advection field, and let
$c_\mu\in L^\infty(\Omega)$ and $f_\mu\in L^2(\Omega)$.  Consider
\begin{equation}\label{eq:transport}
  \divv(\bbeta_\mu u)+c_\mu u=f_\mu\ \text{ in }\Omega,\qquad
  u=g_\mu\ \text{ on }\Gamma_-(\mu):=\{x\in\partial\Omega:\ \bbeta_\mu\cdot\bn(x)<0\}.
\end{equation}
Let $\bSigma_\partial:=\{\btau\in\Hdiv:\ \btau\cdot\bn\in L^2(\partial\Omega)\}$.  We assume
that~\eqref{eq:transport}, with the inflow condition written as
$(\bbeta_\mu u)\cdot\bn=(\bbeta_\mu\cdot\bn)g_\mu$ on $\Gamma_-(\mu)$, has a
unique solution $u(\mu)$ in $\{v\in L^2(\Omega):\ \bbeta_\mu v\in\bSigma_\partial\}$.
Least-squares finite element methods for linear hyperbolic equations were
studied in~\cite{DeSterckEtAl2004} without a flux variable; the
$\Hdiv$-conforming flux was introduced in~\cite{DeSterckEtAl2005} for the
Burgers equation and used for linear transport in~\cite{LiuZhang2020}.  With
$\bsigma=\bbeta_\mu u$, problem~\eqref{eq:transport} becomes
\[
  \bsigma-\bbeta_\mu u=\bzero,\qquad \divv\bsigma+c_\mu u=f_\mu
  \quad\text{in }\Omega,\qquad
  \bsigma\cdot\bn=(\bbeta_\mu\cdot\bn)g_\mu\ \text{ on }\Gamma_-(\mu).
\]
On $\cX_{\rm tr}:=\bSigma_\partial\times L^2(\Omega)$, which does not depend on $\mu$,
define
\begin{equation}\label{eq:Jtr}
  J^{\rm tr}_\mu(\btau,v):=\norm{\btau-\bbeta_\mu v}_0^2
  +\norm{\divv\btau+c_\mu v-f_\mu}_0^2
  +\norm{\btau\cdot\bn-(\bbeta_\mu\cdot\bn)g_\mu}_{0,\Gamma_-(\mu)}^2 ,
\end{equation}
and let $J^{\rm tr,0}_\mu$ be the functional~\eqref{eq:Jtr} with $f_\mu=0$
and $g_\mu=0$.

\begin{prop}\label{prop:transport}
Let $U(\mu):=(\bsigma(\mu),u(\mu))$.  For every $W\in\cX_{\rm tr}$,
\begin{equation}
  J^{\rm tr}_\mu(W)=J^{\rm tr,0}_\mu\bigl(U(\mu)-W\bigr),
\end{equation}
and $J^{\rm tr,0}_\mu(W)=0$ only if $W=0$.  Consequently, for every
finite-dimensional subspace $X_N\subset\cX_{\rm tr}$, the minimizer
$U_N(\mu)$ of $J^{\rm tr}_\mu$ over $X_N$ exists, is unique, and is the best
approximation of $U(\mu)$ in $J^{\rm tr,0}_\mu(\cdot)^{1/2}$, and $\eta_N(\mu):=J^{\rm tr}_\mu(U_N(\mu))^{1/2}
=J^{\rm tr,0}_\mu(U(\mu)-U_N(\mu))^{1/2}$.
\end{prop}

\begin{proof}
The three residuals in~\eqref{eq:Jtr} are affine in $W$ and vanish at
$U(\mu)$.  If $J^{\rm tr,0}_\mu(\btau,v)=0$, then $\btau=\bbeta_\mu v$, so
$\bbeta_\mu v\in\bSigma_\partial$, $\divv(\bbeta_\mu v)+c_\mu v=0$, and
$(\bbeta_\mu v)\cdot\bn=0$ on $\Gamma_-(\mu)$.  By uniqueness, $v=0$ and
$\btau=\bzero$.  Hence $J^{\rm tr,0}_\mu$ is positive definite on $X_N$, and
the rest follows from the projection argument in the proof of
Proposition~\ref{prop:best}.
\end{proof}

Each of the three residuals of a reduced pair is bounded by $\eta_N(\mu)$.
In particular, the balance law for the reduced flux is certified.

\begin{rem}
In the flux formulation with the inflow condition imposed strongly, the
trial space depends on $\mu$.  The span of snapshots computed at parameters
with different inflow boundaries is in general not contained in any of
these spaces.  For example, let $\Omega=(0,1)^2$ and
$\bbeta_\mu=(\cos\mu,\sin\mu)$.  The side $\{x=0\}$ is an inflow side for
$\mu<\pi/2$ and an outflow side for $\mu>\pi/2$, and in the latter case a
snapshot has $\bsigma\cdot\bn=-\cos\mu\,u(\mu)$ on this side.
With~\eqref{eq:Jtr}, the trial space $\cX_{\rm tr}$ is fixed, and the inflow
condition is part of the functional.  For a polyhedral domain with sides
$E$, $\Gamma_-(\mu)$ is the union of the sides with $\bbeta_\mu\cdot\bn_E<0$, and
\[
  \norm{\btau\cdot\bn-(\bbeta_\mu\cdot\bn)g_\mu}_{0,\Gamma_-(\mu)}^2
  =\sum_E\chi_E(\mu)\,
  \norm{\btau\cdot\bn_E-(\bbeta_\mu\cdot\bn_E)g_\mu}_{0,E}^2,
\]
where $\chi_E(\mu)=1$ if $\bbeta_\mu\cdot\bn_E<0$ and $\chi_E(\mu)=0$
otherwise.  This term is parameter-separable with the coefficients
$\chi_E$, $\chi_E\,\bbeta_\mu\cdot\bn_E$, and $\chi_E(\bbeta_\mu\cdot\bn_E)^2$
multiplied by those of $g_\mu$, so the decomposition of
section~\ref{sec:offline} can be used.
\end{rem}

\begin{rem}\label{rem:trnorm}
The certificate is the error in $J^{\rm tr,0}_\mu(\cdot)^{1/2}$, which is
not equivalent to the norm of $\bSigma_\partial\times L^2(\Omega)$~\cite{LiuZhang2020}.
Lemma~\ref{lem:monotone} and Proposition~\ref{prop:sat} only use the
projection property of Proposition~\ref{prop:transport}, so they also hold
for transport.
\end{rem}

\subsection{Snapshots}\label{sec:trsnap}
The scalar variable is in $L^2(\Omega)$, so $\cR_h$ is the identity for
discontinuous Galerkin and finite volume schemes.  For the upwind
discontinuous Galerkin method of degree $k$, the flux does not need a
recovery problem.  Let $\bsigma_h$ in the Raviart--Thomas space $RT_k$
of index $k$ have the normal moments of the
numerical flux on the interior faces, of $(\bbeta_\mu\cdot\bn)g_\mu$ on
$\Gamma_-(\mu)$ (a union of faces), and of $(\bbeta_\mu\cdot\bn)u_h$ on the
other boundary faces, and let its interior moments be those of
$\bbeta_\mu u_h$.  Then $\bsigma_h\in\bSigma_\partial$, and the discontinuous
Galerkin equations give
\begin{equation}
  \divv\bsigma_h+c_\mu u_h-f_\mu=(I-\Pi_k)(c_\mu u_h-f_\mu),
  \qquad
  \bsigma_h\cdot\bn=\Pi_k^\partial\bigl((\bbeta_\mu\cdot\bn)g_\mu\bigr)
  \ \text{ on }\Gamma_-(\mu),
\end{equation}
where $\Pi_k$ and $\Pi_k^\partial$ are the $L^2$ projections onto the piecewise
polynomials of degree $k$ on $\cT_h$ and on the faces of $\Gamma_-(\mu)$,
respectively.  Reconstructions of this type are used
in~\cite{BeckerCapatinaLuce2013,ErnStephansenVohralik2010}.  The snapshot
certificate is computed in one pass over the faces and the elements, and no
linear system is solved.  For piecewise constant $c_\mu$ and data of degree
$k$, only the term $\norm{\bsigma_h-\bbeta_\mu u_h}_0^2$ remains.  For $k=0$,
the method is the cell-centered upwind finite volume scheme, and the face
fluxes of $\bsigma_h\in RT_0$ are its numerical fluxes.  The least-squares
recovery in $RT_k$
gives a certificate that is not larger, at the cost of a symmetric positive
definite problem in $\Hdiv$, and it can also be used for black boxes without
a numerical flux.  At non-intrusive level~0, the truth solver is the least-squares method
of~\cite{LiuZhang2020}.

\subsection{Relation to the double greedy method}

Problem~\eqref{eq:transport} with $\bbeta_\mu=(\cos\mu,\sin\mu)$ is the
transport example of~\cite{DahmenPleskenWelper2014}.  There it is solved in
a variant of the ultraweak formulation of~\cite{DahmenHuangSchwabWelper2012}, in which
the dual norm of the residual equals the $L^2$ error, and the double greedy
algorithm builds a reduced test space together with the trial space to
control this error.  The method of this section certifies the norm of the
residual, that is, the error in $J^{\rm tr,0}_\mu(\cdot)^{1/2}$
(Remark~\ref{rem:trnorm}), with a single greedy loop driven by $\eta_N$.

\section{Numerical experiments}\label{sec:numerics}

\subsection{Setting}\label{sec:setting}
The second-order test problem is the convection--diffusion problem
of~\cite{WuZhang2026} with a two-parameter diffusion coefficient as
in~\cite[Sect.~7]{Zhang2020} and a negative reaction coefficient:
\begin{equation}\label{eq:cdtest}
  -\divv(\nu_\mu\grad u)+\bb\cdot\grad u-3u=1\quad\text{in }\Omega=(0,1)^2,\qquad
  u=0\quad\text{on }\partial\Omega,
\end{equation}
where $\bb=(2,3)^T$, $\nu_\mu=10^{\mu_1}$ for $x<\frac12$, $\nu_\mu=10^{\mu_2}$
for $x>\frac12$, and $\mu\in\cD=[-2,2]^2$.  Thus $10^{-2}\le\nu_\mu\le10^2$, and
$\cJ_\mu(0)=\norm{1}_0^2=1$ for all $\mu$.  The problem has no closed-form
solution, but the certificate does not need one.  Since the reaction
coefficient is negative, the problem is not coercive for small $\nu_\mu$, and
the primal--dual method of~\cite{Zhang2020} cannot be used.
Problem~\eqref{eq:cdtest} is uniquely solvable for all $\mu\in\cD$.  Since $\nu_\mu$ depends continuously on $\mu$
and $\cD$ is compact, the inf--sup constant of $a_\mu$ is bounded below
uniformly on $\cD$.  We use a training
set $\Xi_\train$ of $10^5$ uniformly random points in $\cD$, the first
parameter $(0,0)$, $N_{\max}=40$, and $\rrf=2$.  In
Algorithm~\ref{alg:adaptive}, we use $\varepsilon_h=0.015$
for~\eqref{eq:cdtest} and $\varepsilon_h=0.01$ for the Stokes problem, the
initial mesh
of $8\times8$ squares with each square cut along a diagonal, D\"orfler
marking with $\theta=0.3$ based on the elementwise contributions of
$\cJ_\mu$, and newest vertex bisection.  The space $V_h$ is the continuous
piecewise linear space, and $\bSigma_h$ is the lowest-order Raviart--Thomas
space.  The ``test'' maximum is taken over an independent set
of $10^4$ uniformly random parameters, and the first ten of them are the
test parameters used below.  The black boxes are conforming $P_1$ Galerkin
(non-intrusive level~1), Crouzeix--Raviart with the modified Oswald interpolation
of~\cite{CaiYeZhang2011} (non-intrusive level~2), and finite differences on a uniform grid
(non-intrusive level~3).  The Gram--Schmidt process of section~\ref{sec:offline} uses the
inner products of $\Hdiv$ and $H^1(\Omega)$, row-wise for the Stokes problem,
and of $\Hdiv$ and $L^2(\Omega)$ for transport.  The settings of the Stokes and
the transport problems are given with their results.

In every run below, the maximum of $\eta_N$ over the training set is
nonincreasing, and the saturation trick selects the same parameters as the
full search.  The inequality $\eta_N(\mu_i)\le\eta^\snap_h(\mu_i)$ at the
selected
parameters and, in the adaptive runs, the split~\eqref{eq:split} at ten
test parameters hold up to rounding errors, with a relative error of $2\cdot10^{-6}$ or less.

\Needspace{8\baselineskip}
\subsection{Second-order elliptic equations}

\begin{figure}[htbp]
\centering
\includegraphics[width=\textwidth]{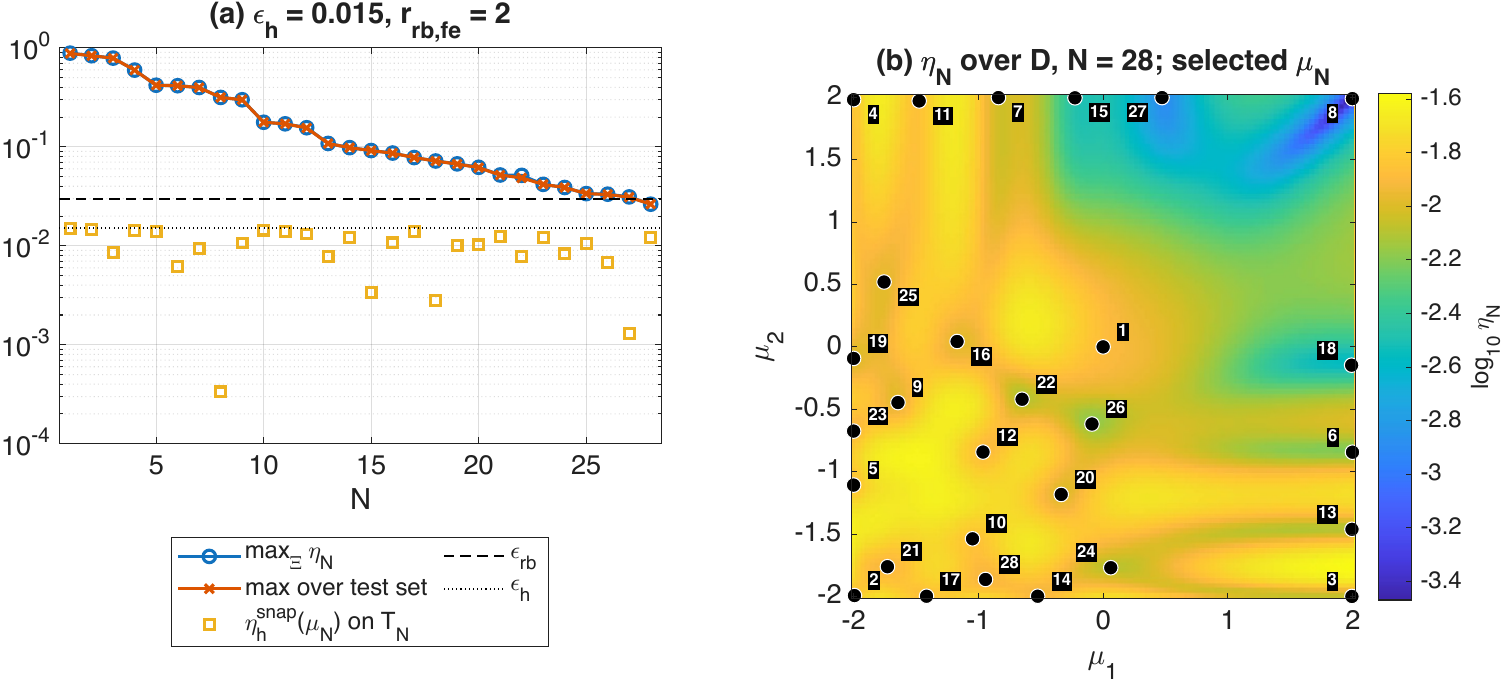}
\caption{Algorithm~\ref{alg:adaptive} at non-intrusive level~0 for~\eqref{eq:cdtest} with
$\varepsilon_h=0.015$ and $\rrf=2$.  (a) The maxima of $\eta_N$ over
$\Xi_\train$ and over the test set, the snapshot certificates
$\eta^\snap_h(\mu_N)=\eta_h(\mu_N)$ on $\cT_N$, $\varepsilon_\rb$, and
$\varepsilon_h$.  (b) $\eta_N$ of the final reduced model over $\cD$ and the
selected parameters (numbered).}
\label{fig:E2}
\end{figure}

\paragraph{Non-intrusive level~0}
The algorithm stops at $N=28$ with
$\max_{\Xi_\train}\eta_{28}=0.0264\le\varepsilon_\rb=0.03$
(Figure~\ref{fig:E2}), and
the maximum over the test set is the same to three digits.  Five of the $28$
selected parameters need mesh refinement.  Most of the refinement occurs at
$\mu_2=(-1.993,-1.993)$, where $\nu_\mu\approx0.01$ on both halves ($38$
refinements, to $376{,}662$ elements), and at $\mu_4=(-2.000,1.981)$ ($20$
more refinements, to $415{,}072$ elements).  The final mesh has $444{,}171$
elements.  The run needs $105$ truth solves, and the
saturation trick evaluates $16\%$ of the certificates of the full search.
The mesh is refined along the interface and along the outflow boundaries
$\{x=1\}$ and $\{y=1\}$.  When $\nu_\mu$ is small on
the left half and large on the right half, the interface is part of the
outflow boundary of the convection-dominated left half, and the solution has
an exponential layer of width about $\nu_\mu/b_1$ on its left.
At ten test parameters, we compute the least-squares solution $U_h(\mu)$ on
the final mesh to obtain the two parts of~\eqref{eq:split}.  The truth part
$\eta_h(\mu)$ is between $4.1\cdot10^{-4}$ and $1.3\cdot10^{-2}$, and the
reduced part $\tnorm{U_h(\mu)-U_N(\mu)}_\mu$ is between $6.3\cdot10^{-4}$
and $1.8\cdot10^{-2}$.  Thus the two parts have comparable ranges.  On the
training set, the stopping test bounds their sum,
$\eta_h^2+\tnorm{U_h-U_N}_\mu^2\le\varepsilon_\rb^2$.

\begin{figure}[htbp]
\centering
\includegraphics[width=\textwidth]{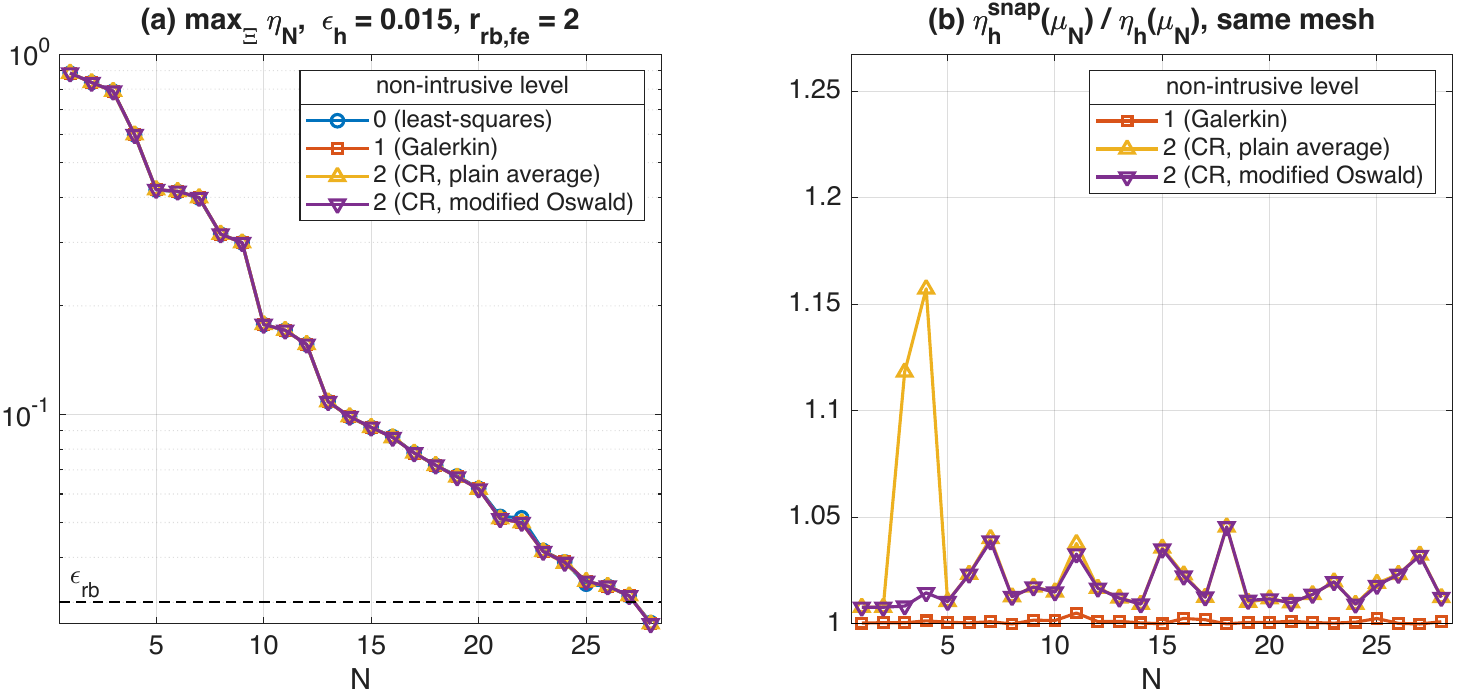}
\caption{Algorithm~\ref{alg:adaptive} for~\eqref{eq:cdtest} with
$\varepsilon_h=0.015$ and $\rrf=2$: non-intrusive level~0 (the run of
Figure~\ref{fig:E2}), non-intrusive level~1 (conforming $P_1$ Galerkin), and non-intrusive level~2
(Crouzeix--Raviart, reconstructed by the plain nodal average and by the
modified Oswald interpolation).  At non-intrusive levels~1 and~2, the flux is given by the
recovery~\eqref{eq:recovery}.  (a) $\max_{\Xi_\train}\eta_N$ and
$\varepsilon_\rb$.  (b) Non-intrusive levels~1 and~2: the snapshot certificate
$\eta^\snap_h(\mu_N)$ divided by the certificate $\eta_h(\mu_N)$ of the
least-squares solution on the same mesh $\cT_N$.}
\label{fig:E3}
\end{figure}

\paragraph{Non-intrusive levels~1--3}
For conforming $P_1$
Galerkin snapshots, the snapshot certificate is at most $0.5\%$ larger than
that of the least-squares solution on the same mesh
(Figure~\ref{fig:E3}(b)).  On a
uniform mesh of $2{,}048$ elements, where the layers are not resolved, this
ratio is $1.56$ at $\mu=(-1.6,0.8)$.  By Proposition~\ref{prop:recbest},
the recovered flux is the best flux for the Galerkin potential, so the
difference only comes from the potential (Remark~\ref{rem:primaryfloor}).
Non-intrusive level~1 stops at $N=28$ with $\max_{\Xi_\train}\eta_{28}=0.0262$, after
$111$ truth solves, on $457{,}760$ elements.  For Crouzeix--Raviart
snapshots with the modified Oswald interpolation, the ratio is between
$1.008$ and $1.046$, and the run stops at $N=28$ with $0.0261$, after $115$
truth solves, on $472{,}466$ elements.  The plain nodal average performs
worse at
$\mu_3\approx(2,-2)$ and $\mu_4\approx(-2,2)$, where $\nu_\mu$ jumps by a
factor of about $10^4$.  Its certificates are $12\%$ and $16\%$ larger than
the least-squares ones, and $\mu_4$ needs $74$ refinements instead of $24$
with the modified Oswald interpolation.
At a vertex on the interface, the plain average also uses values from the side
of the small coefficient, and the flux term of~\eqref{eq:Jelliptic} weights
the resulting gradient error by $\nu_\mu^{1/2}$ on the side of the large
coefficient.  Mesh refinement compensates for this defect, and the final maxima of
$\eta_N$ of the two reconstructions differ by $0.2\%$.  A reduced model
whose snapshots come from non-intrusive levels~0, 1, and~2 in turn also stops at $N=28$,
with $\max_{\Xi_\train}\eta_{28}=0.0258$, after $107$ truth solves.

At non-intrusive level~3, the black box is a finite difference code on its own grid, and
its mesh cannot be adapted.  We therefore run Algorithm~\ref{alg:fixed} on
one uniform grid of $512\times512$ squares.  The solver is the conservative
five-point scheme with central differences for convection (cell P\'eclet
number at most $0.29$), and $\cR_h$ is the $P_1$ interpolation of the grid
values on the grid cut into triangles.  The algorithm starts with
$\varepsilon^0_\rb=10^{-3}$ and raises the tolerance to $\rrf$ times the
largest snapshot certificate.  With finite difference snapshots and with
$P_1$ Galerkin snapshots on the same mesh, the two runs select the same ten
parameters and stop at $N=10$ with $\max_{\Xi_\train}\eta_{10}=0.212$.  The
worst-resolved parameter, $\mu_2=(-1.99,-1.99)$ with
$\eta^\snap_h=0.124$, determines the tolerance for all of $\cD$.

\subsection{The Stokes equations}\label{sec:numstokes}
We solve problem~\eqref{eq:oseen} on $\Omega=(0,1)^2$ with
$\bff_\mu=(-(y-\frac12),x-\frac12)$ and $\nu_\mu=10^{\mu_q}$ on the $q$th of
the four squares of side $\frac12$, where $\mu\in\cD=[-2,2]^4$.  Then
$\tnorm{U(\mu)}_\mu=\cJ_\mu(0)^{1/2}=(\sum_q10^{-\mu_q}/24)^{1/2}$ ranges
from $0.041$ to $4.1$, so we use the relative certificate of
Remark~\ref{rem:relcert}.  The relative certificate does not change when
$\nu_\mu$ is multiplied
by a constant, that is, along $(1,1,1,1)$ in $\cD$, so $\cD$ has three
effective dimensions.  The other settings are those of
section~\ref{sec:setting}.

\begin{figure}[htbp]
\centering
\includegraphics[width=\textwidth]{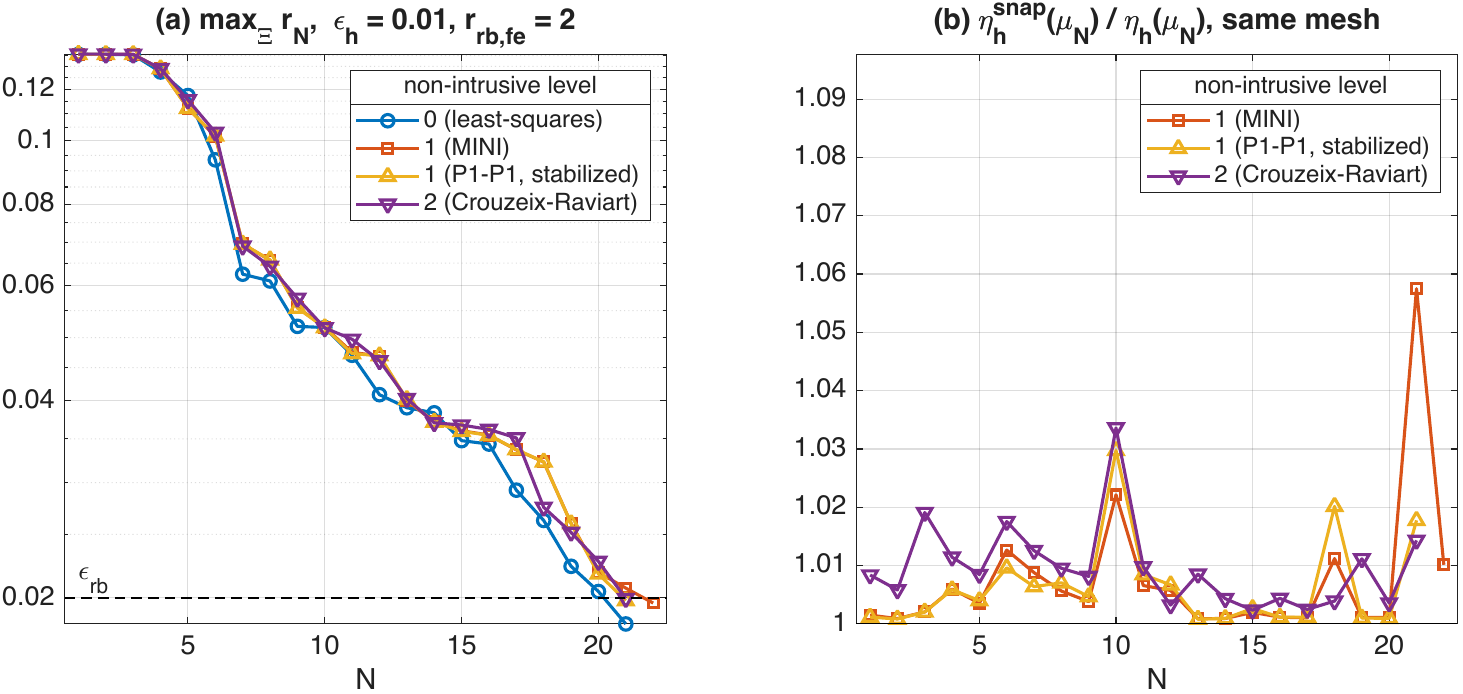}
\caption{The Stokes problem of section~\ref{sec:numstokes}:
Algorithm~\ref{alg:adaptive} with the relative certificate,
$\varepsilon_h=0.01$, and $\rrf=2$, with least-squares snapshots (non-intrusive level~0)
and with the velocity of a black box and the
recovery~\eqref{eq:oseenrec}: MINI and the stabilized $P_1$--$P_1$ element
(non-intrusive level~1), and Crouzeix--Raviart with a jump penalty and the modified Oswald
interpolation (non-intrusive level~2).  The panels are as in Figure~\ref{fig:E3}; in (a),
the certificates are relative.}
\label{fig:E13lev}
\end{figure}

\paragraph{Non-intrusive levels~0--2}
At non-intrusive level~0, the run stops at $N=21$ with
$\max_{\Xi_\train}r_{21}=0.0182$ ($0.0178$ on the test set), on
$22{,}271$ elements, after $49$ truth solves (Figure~\ref{fig:E13lev}).  At
non-intrusive levels~1 and~2, we use three black boxes and take only the velocity from
each.  At non-intrusive level~1, we use the MINI element and the equal-order $P_1$--$P_1$
element.  The latter is not inf--sup stable; it is used with the
Brezzi--Pitk\"aranta stabilization, which perturbs the discrete divergence
constraint.
At non-intrusive level~2, we use the Crouzeix--Raviart velocity with a $P_0$ pressure and
a penalty on the velocity jumps, which is needed for a discrete Korn
inequality for $\beps$, and $\cR_h$ is the modified Oswald interpolation.
The three runs stop at $N=22$, $21$, and $21$, with
$\max_{\Xi_\train}r_N=0.0196$, $0.0198$, and $0.0199$, on final meshes of
$22{,}388$, $21{,}682$, and $23{,}432$ elements.  The ratios of the snapshot
certificate to the non-intrusive level~0 certificate on the same mesh have medians
$1.004$, $1.004$, and $1.008$ and maxima $1.057$, $1.030$, and $1.034$,
respectively (Figure~\ref{fig:E13lev}(b)).

\begin{figure}[htbp]
\centering
\includegraphics[width=0.5\textwidth]{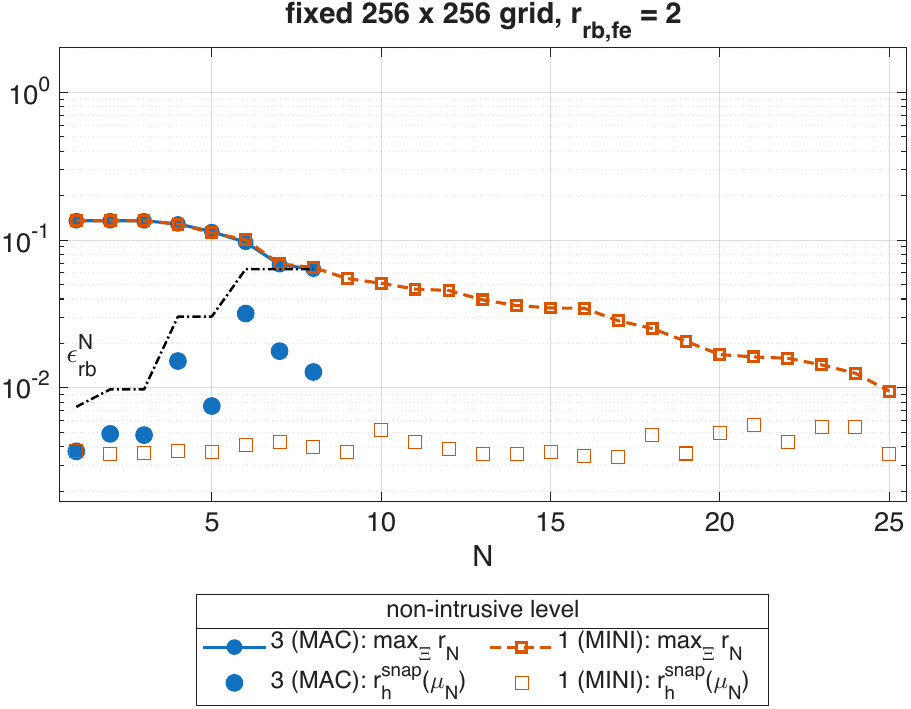}
\caption{The problem of Figure~\ref{fig:E13lev} on one fixed
$256\times256$ grid: Algorithm~\ref{alg:fixed} with relative certificates,
$\varepsilon^0_\rb=10^{-3}$, and $\rrf=2$, with MAC (non-intrusive level~3) and MINI
(non-intrusive level~1) snapshots.  Shown are $\max_{\Xi_\train}r_N$, the relative
snapshot certificates
$r^\snap_h(\mu_N):=\eta^\snap_h(\mu_N)/\cJ_{\mu_N}(0)^{1/2}$, and
$\varepsilon^N_\rb$ of the MAC run.}
\label{fig:E13fixed}
\end{figure}

\paragraph{Non-intrusive level~3}
Figure~\ref{fig:E13fixed} compares MAC and MINI snapshots on one fixed
$256\times256$ grid.  We use the MAC scheme in finite volume form, with
the symmetric gradient, $\nu_\mu$ at the cell centers, and the
mean of the neighboring cells at the vertices; $\cR_h$ is as in
section~\ref{sec:oseen}.  The relative snapshot certificates of the MAC run are
between $0.0037$ and $0.032$.  The largest one is at $\mu_6$, where MINI on
the same grid gives $0.0041$.  This snapshot raises $\varepsilon^N_\rb$ to
$0.064$, and
the run stops at $N=8$ with $\max_{\Xi_\train}r_8=0.064$.  Up to $N=8$, the
maxima of $r_N$ of the MAC and MINI runs are nearly the same.  With MINI
snapshots, the relative snapshot certificates are between $0.0034$ and
$0.0056$, and the run stops at $N=25$ with $0.0095$.  For both solvers, the
tolerance is determined by the worst-resolved selected parameter.

\Needspace{8\baselineskip}
\subsection{Linear transport equations}

\begin{figure}[htbp]
\centering
\includegraphics[width=\textwidth]{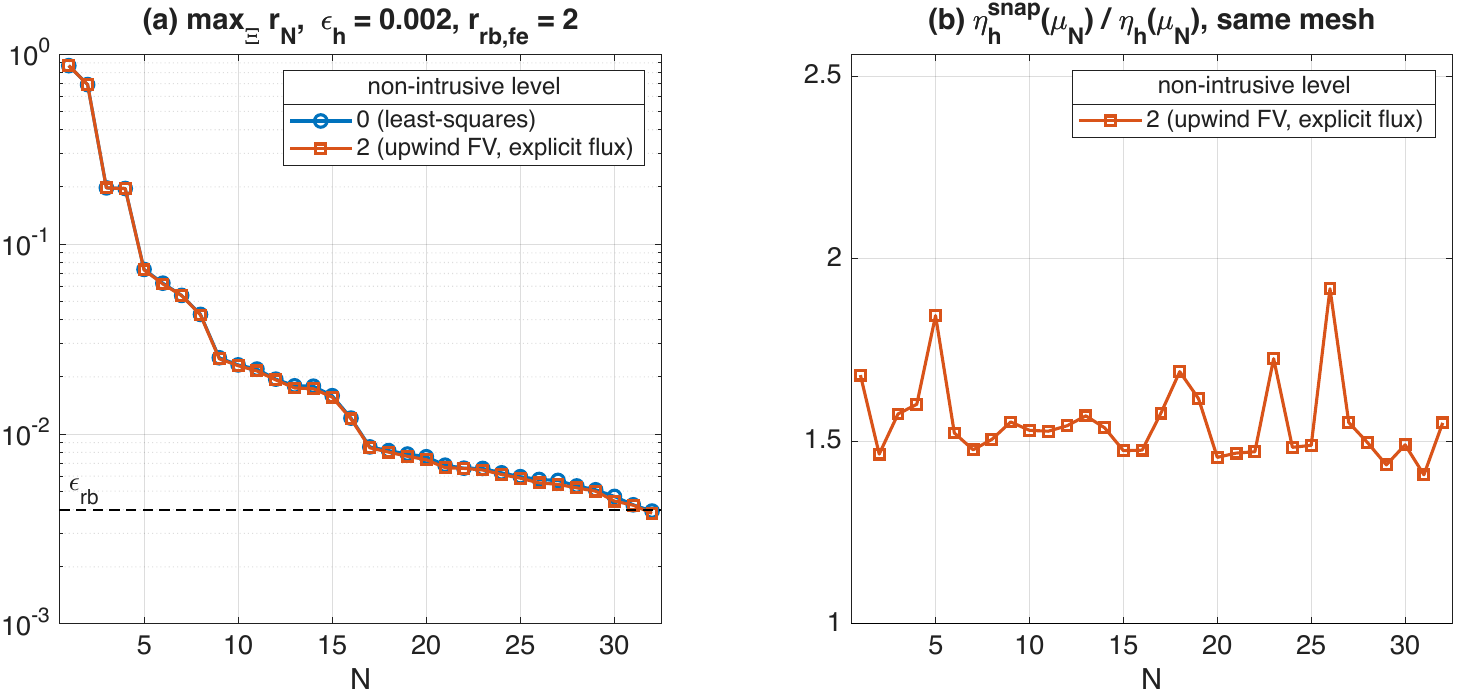}
\caption{Problem~\eqref{eq:transport} with
$\bbeta_\mu=(\cos\mu,\sin\mu)$, $\mu\in[0,\pi]$, and the data
of~\cite{DahmenPleskenWelper2014}: Algorithm~\ref{alg:adaptive} with the relative
certificate, $\varepsilon_h=2\cdot10^{-3}$, and $\rrf=2$, with
least-squares snapshots (non-intrusive level~0) and with upwind finite volume snapshots
and the flux of section~\ref{sec:trsnap} (non-intrusive level~2).  The panels are as in
Figure~\ref{fig:E3}; in (a), the certificates are relative.}
\label{fig:E15lev}
\end{figure}

\begin{figure}[htbp]
\centering
\begin{minipage}[b]{0.48\textwidth}\centering
\includegraphics[width=\textwidth]{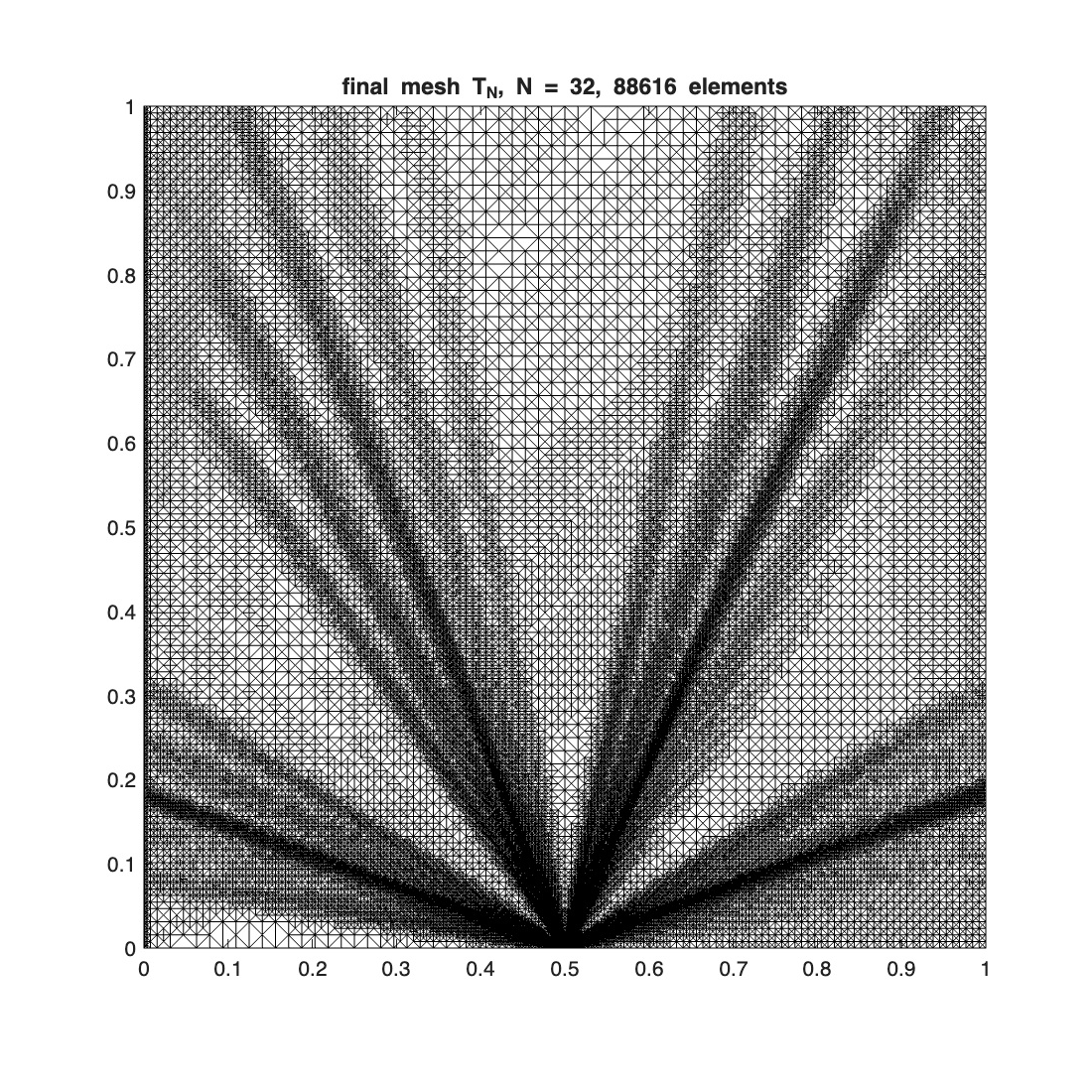}\\(a)
\end{minipage}\hfill
\begin{minipage}[b]{0.48\textwidth}\centering
\includegraphics[width=\textwidth]{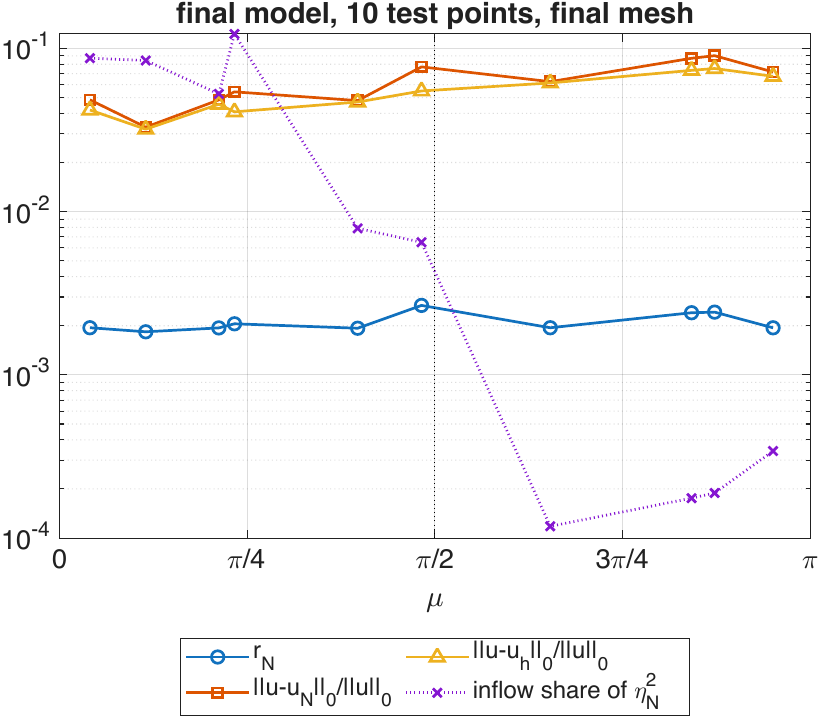}\\(b)
\end{minipage}
\caption{The non-intrusive level~0 run of Figure~\ref{fig:E15lev}.  (a) The final mesh
with $88{,}616$ elements.  (b) At ten test parameters on the final mesh:
the relative certificate $r_{32}$ of the final model, the $L^2$ errors of
$u_{32}$ and of the least-squares solution $u_h$ divided by $\norm{u}_0$,
and the share of the inflow term in $\eta_{32}^2$.}
\label{fig:E15mesh}
\end{figure}

\paragraph{Non-intrusive level~0}
We solve~\eqref{eq:transport} on $\Omega=(0,1)^2$ with
$\bbeta_\mu=(\cos\mu,\sin\mu)$, $\mu\in\cD=[0,\pi]$, and the data
of~\cite{DahmenPleskenWelper2014}: $c=1$, $f=\frac12$ for
$x<y$ and $f=1$ elsewhere, and $g=1-y$ for $x\le\frac12$ and $g=0$
elsewhere.  The solution $u(\mu)$ jumps across the characteristic through
$(\frac12,0)$, which turns with $\mu$.  Since this line is parallel to
$\bbeta_\mu$, the normal component of $\bsigma=\bbeta_\mu u$ is continuous
across it, and $\bsigma\in\bSigma_\partial$.  The side $\{x=0\}$ is an inflow side
for $\mu<\pi/2$, and the side $\{x=1\}$ is an inflow side for
$\mu>\pi/2$.  Non-intrusive level~0 uses $RT_0\times P_0$, which contains the non-intrusive level~2
pairs.  The line $\{x=y\}$ consists of edges of the initial mesh of
$8\times8$ squares, and $(\frac12,0)$ is a vertex of this mesh, so $J^{\rm tr}_\mu(0)$ is integrated exactly on
every mesh.  Over the selected parameters, $J^{\rm tr}_\mu(0)^{1/2}$ is between
$0.79$ and $1.06$.  We use $\varepsilon_h=2\cdot10^{-3}$, the first
parameter $\pi/2$, and training and test sets of $10^4$ and $10^3$
uniformly random points.  The run stops at $N=32$ after $184$ truth solves, on
$88{,}616$ elements, with $\max r_{32}=3.93\cdot10^{-3}$ on both the
training and the test set (Figure~\ref{fig:E15lev}).  The first three
parameters are $\pi/2$ and the training points closest to $0$ and $\pi$.  The mesh is refined along the
characteristics through $(\frac12,0)$ at the selected parameters, which form
a fan of $32$ lines (Figure~\ref{fig:E15mesh}(a)).

The $L^2$ errors with respect to the exact solution are not certified
(Remark~\ref{rem:trnorm}).  The exact solution is evaluated along the
characteristics, and the $L^2$ errors are computed by a composite rule with
$16$ subtriangles in each element and a $12$-point rule on each.  At ten test parameters, $r_{32}$ is between
$1.8\cdot10^{-3}$ and $2.7\cdot10^{-3}$.  For $U_{32}=(\bsigma_{32},u_{32})$,
$\norm{u-u_{32}}_0/\norm{u}_0$ is between $0.033$ and $0.090$, and for the
least-squares solution on the final mesh, $\norm{u-u_h}_0/\norm{u}_0$ is
between $0.032$ and $0.075$ (Figure~\ref{fig:E15mesh}(b)).  The share of the inflow term in
$\eta_{32}^2$ is between $0.05$ and $0.12$ at the test parameters below
$\pi/4$ and between $1.2\cdot10^{-4}$ and $3.4\cdot10^{-4}$ at those above
$\pi/2$, where the inflow side $\{x=1\}$ carries the datum $g=0$.  The $L^2$
error of $u_{32}$ is at most
$1.41$ times that of $u_h$ on the same mesh, so the gap between the $L^2$
error and the certificate is already present in the truth solution: the
certificate is the error in $J^{\rm tr,0}_\mu(\cdot)^{1/2}$, which is not
equivalent to the $L^2$ norm (Remark~\ref{rem:trnorm}).  On uniform meshes of $16\times16$ to
$64\times64$ squares, at $\mu=0.6$ and $\mu=2.2$, $\eta_h$ decreases at the
rates $0.92$ and $0.73$ in $h$, and the relative $L^2$ error at the rates
$0.36$ and $0.30$.  For $f=1$ and $g=0$, which is also an
example of~\cite{DahmenPleskenWelper2014} and where $u$ only has a kink, the
rates are $1.00$ and $0.98$ for $\eta_h$ and $0.99$ and $0.89$ for the $L^2$
error.

\paragraph{Non-intrusive level~2}
We use the upwind finite volume scheme, that is, the case
$k=0$ of section~\ref{sec:trsnap}.  Its numerical fluxes give the $RT_0$
flux directly, without a linear solve.  The run stops at $N=32$, as at
non-intrusive level~0, with $\max_{\Xi_\train}r_{32}=3.83\cdot10^{-3}$
(Figure~\ref{fig:E15lev}).  The two levels select the same parameters to
within $3\cdot10^{-3}$, except that steps~14 and~15 and steps~27 and~28 are
interchanged.  The finite volume snapshots are less accurate: their
certificate is $1.41$ to $1.92$ times (median $1.53$) the non-intrusive level~0
certificate on the same mesh.  Hence the final mesh has $174{,}447$ elements
instead of $88{,}616$, and the run needs $271$ truth solves instead of
$184$.  For
$f=1$ and $g=0$, with $\varepsilon_h=10^{-3}$, both levels stop at $N=20$,
the median of the ratio is $1.67$, and the final meshes have $63{,}354$ and
$157{,}679$ elements.

\section{Concluding remarks}
For a first-order system least-squares formulation, the square root of
the functional at any conforming pair is the norm of the residual of the
continuous problem at the pair.  This separates the computation of the snapshots
from the reduction and the certification.  The truth solver only affects
the accuracy of the snapshots, which their certificates measure; the
online reduced problem, the certificate, and the greedy algorithms are the same
for all solvers.  The online reduced problem needs no inf--sup condition on
the reduced spaces, and its error certificate needs no stability constant.  In the experiments with adaptive meshes, the reduced dimension hardly
depends on the solver, while a less accurate solver needs finer truth
meshes.

\end{document}